\documentclass[11pt,leqno]{article}
\usepackage[T1]{fontenc}
\usepackage{lmodern,amsmath,amssymb,amsthm,amscd,mathtools}
\usepackage[margin=28mm]{geometry}
\usepackage{microtype,longtable,booktabs,array}
\usepackage[colorlinks=true,linkcolor=blue,citecolor=blue,urlcolor=blue]{hyperref}
\usepackage[utf8]{inputenc}
\numberwithin{equation}{section}
\newtheorem{theorem}{Theorem}[section]
\newtheorem{proposition}[theorem]{Proposition}
\newtheorem{lemma}[theorem]{Lemma}
\newtheorem{corollary}[theorem]{Corollary}
\theoremstyle{definition}\newtheorem{definition}[theorem]{Definition}
\theoremstyle{remark}\newtheorem{remark}[theorem]{Remark}
\theoremstyle{definition}
\newtheorem{example}[theorem]{Example}
\newcommand{\C}{\mathbb C}\newcommand{\Q}{\mathbb Q}\newcommand{\Z}{\mathbb Z}
\newcommand{\HH}{\mathbb H}\newcommand{\T}{\mathsf T}\newcommand{\G}{\mathsf G}
\newcommand{\R}{\mathsf R}\newcommand{\A}{\mathsf A}
\newcommand{\dd}{\mathrm d}\newcommand{\dR}{\mathrm{dR}}
\newcommand{\Spec}{\operatorname{Spec}}\newcommand{\Res}{\operatorname{Res}}
\newcommand{\Reg}{\operatorname{Reg}}\newcommand{\Gr}{\operatorname{Gr}}
\newcommand{\im}{\operatorname{im}}
\newcommand{\tr}{\mathsf t}
\newcommand{\cO}{\mathcal O}

\newcommand{\Gm}{\mathbb G_m}\newcommand{\Ga}{\mathbb G_a}
\newcommand{\SL}{\operatorname{SL}}\newcommand{\Sp}{\operatorname{Sp}}
\newcommand{\wt}{\operatorname{wt}}\newcommand{\ind}{\operatorname{ind}}
\newcommand{\eps}{\varepsilon}

\newcommand{\subjclass}[2][2020]{%
  \def\GMCDmsc{\textit{#1 Mathematics Subject Classification.}\enspace #2.}%
  \hypersetup{pdfsubject={MSC #1: #2}}}
\newcommand{\keywords}[1]{%
  \def\GMCDkeywords{\textit{Key words and phrases.}\enspace #1.}%
  \hypersetup{pdfkeywords={#1}}}


\title{Gauss-Manin connection in disguise:\\
limiting periods and quasi Jacobi forms}
\author{Jin Cao \thanks{School of Mathematics and Physics, University of Science and Technology Beijing, Beijing, China \\
\textit{E-mail address}: \texttt{caojin@ustb.edu.cn}
  \newline\GMCDmsc\newline\GMCDkeywords}}
\date{}

\begin{document}
\maketitle
\begin{abstract}
We construct an algebraic moduli space of dimension six parametrizing elliptic curves \(X\) with origin \(O\), a marked point \(P\ne O\), a relative de Rham frame, and a nonzero tensor in \(T_OX\otimes T_PX\). Its normalized period map is a biholomorphism. We evaluate its inverse on a family of period matrices parametrized by \((\tau,z)\), where \(z\) represents \(P\) on \(\C/(\Z\tau+\Z)\). The tensor coordinate is then \((\theta_1/\eta^3)^2\); the other coordinates are given by Weierstrass functions, Eisenstein series, and a logarithmic theta derivative. Two Gauss-Manin vector fields realize differentiation in \(\tau\) and \(z\) on quasi Jacobi forms of integral index. For an explicit genus two smoothing, regularized periods identify the punctured normal bundle of the nodal divisor with the bundle of nonzero tangent tensors at the two points over the node. We compute the limiting connection, the boundary values of the genus two generators, and the unit relating the smoothing and Fourier-Jacobi parameters.

\end{abstract}

\section{Introduction}\label{sec:introduction}

\subsection{The comparison problem}

The Gauss-Manin connection in disguise (GMCD) encodes the differential equations satisfied by modular functions as vector fields on moduli spaces enhanced by de Rham frames. The elliptic construction yields the Ramanujan equations \cite{Mov,MovBook}; its relative version for elliptic curves with two marked points yields quasi Jacobi forms of index zero \cite{CMVL}. The genus two construction of \cite{CMY} gives the enhanced moduli space and its differential equations for Siegel modular forms.

We consider a family of smooth genus two curves degenerating to a curve with one nonseparating node. The normalization of the singular fiber is a connected elliptic curve \(X\), and the two points \(O,P\in X\) above the node are identified in the singular curve. Let \(W\) denote the weight filtration on the limiting mixed Hodge structure of the first rational cohomology of the smooth fibers. Its graded pieces in weights \(0,1,2\) are \(\Q(0), H^1(X,\Q), \Q(-1)\), respectively. The weight zero piece comes from the first cohomology of the dual graph. In the logarithmic de Rham description, the weight two piece is detected by the residue on one branch of the node; the residue on the other branch has the opposite sign. The cup product pairs the weight zero and weight two pieces perfectly. To compare the de Rham frames we also need the period along a cycle meeting the vanishing cycle once. After subtraction of its logarithmic term, the finite part depends on the smoothing parameter and the tangents at \(O\) and \(P\).

To recover the limiting genus two periods from the elliptic normalization \(X\), we compare the finite part of the period along the crossing cycle with the integral from \(O\) to \(P\) regularized with respect to the chosen endpoint tangents.  We use the classical descriptions of the limiting mixed Hodge structure, the periods under plumbing, and the Fourier-Jacobi expansion \cite{Steenbrink,Hain,Yamada,HN,EZ,vdG}. For the explicit smoothing in Section~\ref{sec:degeneration-new}, we identify the limiting frame and connection with those of an enhanced elliptic curve.  This comparison determines the unit relating the algebraic smoothing parameter to the Siegel period coordinate. We compute the boundary values of the generators of \cite[(4.11)]{CMY} in the original genus two frame. A single multiplicative coordinate records the dependence on the endpoint tangents and on the coefficient of the graph class in the chosen lift of the residue class.

\subsection{Enhanced moduli spaces}
The space \(\T_0\) parametrizes elliptic curves with two marked points and a relative de Rham frame. In \(\T^{\log}\) one adds a class with residues \(-1,1\) at \(O,P\) that completes it to a symplectic frame; the holomorphic differential fixes the two endpoint tangents. In \(\T\) these tangents may be chosen independently. In \(\T^\times\) the additional datum is a single nonzero tensor in \(T_OX\otimes T_PX\).

Throughout the paper we work over $\C$. Put
\begin{equation}\label{eq:base}
 t_3=4a^3-t_2a-b^2,\qquad
 \Delta=27t_3^2-t_2^3,\qquad
 A_0=\C[a,b,c,t_1,t_2,\Delta^{-1}].
\end{equation}
The elliptic curve \(X\), its origin \(O\), its second marked point \(P\), and the invariant differential \(\alpha_1\)
are given by
\begin{equation}\label{eq:curve}
 X:\quad y^2=4x^3-t_2x-t_3,\qquad
 O=[0:1:0],\qquad P=(a,b),\qquad \alpha_1=\frac{\dd x}{y}.
\end{equation}
The differential $\alpha_1$ determines normalized endpoint tangents
\begin{equation}\label{eq:normalized-tangents-intro}
 \partial_O=(\alpha_1|_O)^\vee,\qquad
 \partial_P=(\alpha_1|_P)^\vee.
\end{equation}
Write \(v_O=\lambda\partial_O\) and \(v_P=\mu\partial_P\), with \(\lambda,\mu\in\C^\times\). Normalize the graph class
\(\alpha_0\) to have period \(-1\) along a path from \(O\) to \(P\). We write the last element of the frame as
\(\alpha_3=\eta_P+c\alpha_1+d_T\alpha_0\), where \(\eta_P\) is the differential of the third kind in
Section~\ref{sec:moduli}, with residues \(-1\) at \(O\) and \(1\) at \(P\). On the section \(\lambda=\mu=1\), write
\(d=d_T\). The four parameter spaces are
\begin{equation}\label{eq:four-enhanced-spaces}
\begin{aligned}
\T_0&=\Spec A_0, &\dim_{\C}\T_0&=5,\\
\T^{\log}&=\Spec A_0[d], &\dim_{\C}\T^{\log}&=6,\\
\T&=\Spec A_0[d_T,\lambda^{\pm1},\mu^{\pm1}], &\dim_{\C}\T&=8,\\
\T^\times&=\Spec A_0[h,h^{-1}], &\dim_{\C}\T^\times&=6.
\end{aligned}
\end{equation}
Their moduli interpretations are as follows.
\begin{enumerate}
\item \(\T_0\) parametrizes the relative de Rham frame \((\alpha_0,\alpha_1,\alpha_2)\) of \cite{CMVL}. The graph class
\(\alpha_0\) can be represented by \(\dd f\) with \(f(O)=1\), \(f(P)=0\).\footnote{Here \(f\) is a rational function
regular at \(O,P\); its other poles are treated in the \v Cech-de Rham complex.} The parameters \(c,t_1\) determine
\(\alpha_2\) by \eqref{eq:companion}.
\item \(\T^{\log}\) adds a class \(\alpha_3\) that completes the relative frame to a symplectic frame. The endpoint
tangents are fixed by \eqref{eq:normalized-tangents-intro}.
\item \(\T\) records this symplectic frame together with two independently chosen nonzero endpoint tangents
(Definition~\ref{def:full-T}).
\item $\T^\times$ parametrizes the relative enhancement together with
a nonzero tensor in $T_OX\otimes T_PX$
(Definition~\ref{def:effective-T}).
\end{enumerate}
Fix \((X,O,P;\alpha_0,\alpha_1,\alpha_2)\), a homology marking, and a path from \(O\) to \(P\). Let \(K\) be the
regularized integral of \(\alpha_3\) along this path at \((d_T,\lambda,\mu)=(0,1,1)\). Continuing from these values,
the integral is \(K-d_T+\log\lambda+\log\mu\), with both logarithms continued from \(\log1=0\). Thus the effective
coordinate is
\begin{equation}\label{eq:effective-intro}
 h=\lambda\mu e^{-d_T}.
\end{equation}

Formula \eqref{eq:effective-intro} defines a holomorphic map \(\T^{\mathrm{an}}\to\T^{\times,\mathrm{an}}\). The
tensor description gives \(\T^\times\) its algebraic structure. On the section \(\lambda=\mu=1\), we have \(d=d_T\)
and \(h=e^{-d}\). Conversely, every point of \(\T^{\times,\mathrm{an}}\) has a simply connected analytic neighborhood
on which one can choose \(d=-\log h\); any two choices differ by a constant in \(2\pi i\Z\).

The group \(\G\) of changes of de Rham frame acts algebraically on \(\T\). If the change of frame multiplies
\(\alpha_1\) by \(s\) and translates \(d_T\) by \(w\), its action on \(\T^\times\) is \(h'=s^2e^{-w}h\). Thus the action
on \(h\) is holomorphic, and with compatible logarithms it gives \(d'=d+w-2\log s\). Writing \(\G=\G_0\times\Ga\), we
obtain an algebraic action on \(\T^\times\) by replacing \(e^{-w}\) with an independent coordinate of \(\Gm\). The
homomorphism \(\G^{\mathrm{an}}\to(\G_0\times\Gm)^{\mathrm{an}}\) has kernel \(\{1\}\times2\pi i\Z\) and identifies the
infinitesimal actions (Proposition~\ref{prop:effective-algebraic-group}). The connection and its modular vector fields
are algebraic on \(\T^\times\).

\subsection{The period map, regular functions, and modular vector fields}
We formulate the period isomorphism, identify the pullbacks of regular
functions under the  \(\mathsf t\)-map, and compare the Gauss-Manin
vector fields with analytic differential operators.

Let \(\Pi\) be the generalized period domain, and let \(\Gamma\) be its group of changes of Betti homology marking, defined in Section~\ref{sec:periods}. Fix \(k=(-2\pi i)^{-1/2}\). We write \(F=\theta_1/\eta^3\) for the odd Jacobi theta function \(\theta_1\) divided by the cube of the Dedekind eta function \(\eta\), with the conventions of Section~\ref{sec:analytic-Jacobi}. For a point of \(\T^{\times,\mathrm{an}}\), choose a homology marking, a lift of \(P\), and a local logarithm \(d=-\log h\). Integrate the corresponding de Rham frame over the marked cycles, regularizing the integral from \(O\) to \(P\) with the tangents in \eqref{eq:normalized-tangents-intro}. Normalize the matrix as in \eqref{eq:Q}; changing these choices leaves its class modulo \(\Gamma\) unchanged.
\begin{theorem}\label{thm:theta}
The resulting period map is a biholomorphism
\[
 \operatorname{Per}:\T^{\times,\mathrm{an}}
       \xrightarrow{\ \sim\ }\Gamma\backslash\Pi.
\]
\end{theorem}

Forgetting the tangent tensor corresponds to deleting the last row and
column of the period matrix. On the base, the resulting period map is the
index zero isomorphism of \cite[\S8.3, Proposition~8.1]{CMVL}. On each
fiber, let \(w_*\) denote the period coordinate at \(h=1\). The map
\(h\mapsto[w-w_*]=[k^2\log h]\) is a biholomorphism from
\(\C^\times\) to \(\C/\Z\) (Theorem~\ref{thm:inverse}). We call the image of \((\tau,z)\mapsto[\mathcal P(\tau,z)]\) in
\(\Gamma\backslash\Pi\) the classical period locus, where \(\mathcal P(\tau,z)\) is the matrix in
\eqref{eq:tau-locus}. Composing this parametrization with \(\operatorname{Per}^{-1}\) gives
\[
\mathcal H\xrightarrow{\mathsf t}\T^{\times,\mathrm{an}},
\qquad
\mathcal H=\{(\tau,z)\in\HH\times\C:z\notin\Z\tau+\Z\}.
\]
Write \(J_1=\partial_z\log F\).  Let \(\wp(\tau,z)\) be the Weierstrass function associated
with the lattice \(\Z\tau+\Z\), and put \(\wp'=\partial_z\wp\). 
Let \(E_2(\tau)\) and \(E_4(\tau)\) denote the normalized Eisenstein series of weights \(2\) and \(4\), respectively, as in Section~\ref{sec:analytic-Jacobi}. 


\begin{theorem}\label{thm:intro-coordinate-pullbacks}
Along the $\mathsf t$-map, one has
\begin{equation}\label{eq:theta-map}
\begin{alignedat}{3}
 \mathsf t^*a&=k^2\wp(\tau,z),
 &\qquad \mathsf t^*b&=k^3\wp'(\tau,z),
 &\qquad \mathsf t^*c&=-kJ_1(\tau,z),\\
 \mathsf t^*t_1&=-\frac{k^{-2}}{12}E_2(\tau),
 & \mathsf t^*t_2&=\frac{k^{-4}}{12}E_4(\tau),
& \qquad  \mathsf t^*h&=F(\tau,z)^2.
\end{alignedat}
\end{equation}
For every \(f\in\C[a,b,c,t_1,t_2,h,\Delta^{-1},h^{-1}]\), substitution from \eqref{eq:theta-map} gives \(\mathsf t^*f\); the pullbacks of \(\Delta\) and \(h\) are nowhere zero on \(\mathcal H\).
\end{theorem}
We prove these formulas using
automorphy, the Gauss-Manin connection, and tangential normalization.  See Proposition~\ref{prop:t-map-coordinates}.
In particular, $\mathsf t^*h$ extends to a weak Jacobi form of weight
$-2$ and index $1$.

The period identification also respects the differential equations.
\begin{theorem}\label{thm:main}
\begin{enumerate}
\setlength{\itemsep}{0pt}
\item The Gauss-Manin connection matrix \(\A\) on \(\T^{\log}\) is algebraic, flat, and preserves \eqref{eq:Phi}.
Substituting \(\dd d=-\dd h/h\) gives a flat algebraic connection on \(\cO_{\T^\times}^{\oplus4}\).
\item On either space, the equations \(\A(\R_i)=C_i\), with \(C_i\) as in \eqref{eq:matrices}, determine a global
frame of six regular vector fields. Their \(\C\) linear span is the Jacobi Lie algebra
\(\mathfrak{sl}_2\ltimes\mathfrak h_3\).
\item Along \(\mathsf t\),
\[
\dd\mathsf t(\partial_\tau)=\R_\tau\circ\mathsf t,\qquad
\dd\mathsf t(k\partial_z)=\R_z\circ\mathsf t.
\]
\end{enumerate}
\end{theorem}
Here \(\mathfrak h_3\) is the three dimensional Heisenberg Lie algebra.
Sections~\ref{sec:connection}--\ref{sec:vector-fields} compute the connection and vector fields;
Corollary~\ref{prop:geometric} and Section~\ref{sec:Jacobi} give their period realization.

\subsection{Degeneration and boundary values}

Choose \(R=(A,B)\in X\) with \(2R=P\). Work on a sufficiently small polydisk \(\mathcal B\) with coordinates \((A,c,t_1,t_2,t_3)\), choosing the local branch of \(B^2=4A^3-t_2A-t_3\) through \(R\) and assuming \(B\Delta\ne0\). Over \(S=\mathcal B\times\mathbb D\), consider the family
\[
 C_s:\quad Y^2=((X-A)^2-s)(4X^3-t_2X-t_3),\qquad
 \rho=-\frac{s}{4B^2}.
\]
Along \(D=\{s=0\}\), the two points \(-R,R\) on the normalization are identified. Translation by \(R\) carries this ordered pair to \((O,P)\).

For the period comparison, we work in the complex analytic topology. Over a simply connected sector of the punctured \(s\) disk, choose a locally constant symplectic basis \((A_1,A_2,B_1,B_2)\) of \(H_1(C_s,\Z)\), with \(A_2\) the vanishing cycle oriented positively about \(P\). Let \(P_2\) be the period matrix of the frame in Theorem~\ref{thm:enhanced-nonsep-limit}, with rows ordered as \((B_2,B_1,A_1,A_2)\), and let \(N\) be the logarithm of positive monodromy on these rows. On the normalization put
\[
 \omega_1=\int_{B_1}\alpha_1,\qquad \omega_2=\int_{A_1}\alpha_1,\qquad
 v=2\int_O^R\alpha_1,\qquad \tau=\frac{\omega_1}{\omega_2},\qquad
 z=\frac{v}{\omega_2},
\]
where the path from \(O\) to \(R\) is fixed. Since \(2R=P\), \(v\) is a lift of \(P\). Choose \(B_2\) so that, after translation by \(R\), its limiting path lifts from \(0\) to \(v\). Use the same elliptic cycle lifts and endpoint logarithms for \(P_2\) and \(Q\), with \(\log\rho\) chosen compatibly with the return of \(B_2\) through the annulus.

The elliptic period matrix \eqref{eq:Q} at \(d=0\) is
\[
 Q(a,b,c,t_1,t_2,0)=
 \begin{pmatrix}
 -1&v&t_1v-g&w_0\\
 0&\omega_1&t_1\omega_1-\eta_1&g\omega_1-\eta_1v\\
 0&\omega_2&t_1\omega_2-\eta_2&g\omega_2-\eta_2v\\
 0&0&0&2\pi i
 \end{pmatrix},
\]
where \(\Lambda=\Z\omega_1+\Z\omega_2\), \(g=c+\zeta_\Lambda(v)\), \(w_0=-2\log\sigma_\Lambda(v)+gv-2\log(2\pi k)\), and \(\eta_j\) is the quasi-period of \(\zeta_\Lambda\) corresponding to \(\omega_j\). Finally, write \(q_2=e^{2\pi i\Omega_{22}}\), where \(\Omega\) is the Siegel period matrix for the chosen marking.
\begin{theorem}\label{thm:intro-degeneration}
For this enhanced family the following hold.
\begin{enumerate}
\setlength{\itemsep}{0pt}
\item \(N^2=0\) and \(\operatorname{rank}N=1\). The regularized period matrix extends holomorphically across \(D\), and its restriction to \(D\) is invertible. More precisely,
\begin{equation}\label{eq:intro-boundary-comparison}
\widehat P_2=e^{-(\log\rho)N/(2\pi i)}P_2,\qquad
\widehat P_2|_D=Q(a,b,c,t_1,t_2,0).
\end{equation}
\item Put \(D_0=\operatorname{diag}(1,k,k,k^2)\). The map \(\ell:D\to\T^{\times,\mathrm{an}}\) defined by
\(\ell(u)=\operatorname{Per}^{-1}([\widehat P_2(u)D_0])\) satisfies
\[\ell^*h=1,\qquad
\left.(\widehat P_2^{-1}\dd\widehat P_2)^{\tr}\right|_D=\ell^*\A.
\]
Allowing the nonzero normal vector to vary gives an isomorphism
\(N_{D/S}^{\times}\simeq\ell_0^*(\T^{\times,\mathrm{an}})\) over \(D\), where \(\ell_0=p\circ\ell\) and
\(p:\T^\times\to\T_0\) forgets the effective tensor.
\item The quotient \(q_2/\rho\) extends to a holomorphic unit with
\begin{equation}\label{eq:intro-normal-unit}
\left.\frac{q_2}{\rho}\right|_D=-\frac{4\pi^2}{\omega_2^2 F(\tau,z)^2}.
\end{equation}
\end{enumerate}
\end{theorem}
The parameter \(\rho\) selects the normal vector \(\xi\) with \(\dd\rho(\xi)=1\) at each point of \(D\).
The associated elliptic enhancement has \(h=1\); the map
\(\mathsf t\) gives \(\mathsf t^*h=F^2\) by Theorem~\ref{thm:intro-coordinate-pullbacks}.
If \(\rho'=u\rho+O(\rho^2)\) with \(u|_D\ne0\), the normal vector selected by \(\rho'\)
is \((u|_D)^{-1}\xi\). Its effective coordinate and the boundary value of \(q_2/\rho'\)
are both scaled by \((u|_D)^{-1}\). Theorem~\ref{thm:normal-bundle} gives the connection
on the full punctured normal bundle.


Proposition~\ref{prop:CMY-explicit} computes the boundary values of the genus two generators
in the frame \eqref{eq:CMY-full-frame}; Corollary~\ref{cor:invariant-boundary-values} identifies combinations that
are independent of the choice of \(R\). To express a leading algebraic term in period coordinates,
one uses the unit \eqref{eq:intro-normal-unit} and, for a scalar modular form, the change of
frame \eqref{eq:algebraic-classical-frame-conversion}. At a separating node the monodromy
on \(H^1\) is trivial. The limiting connection is the direct sum of the two elliptic connections;
with a split frame the boundary period matrix is block diagonal
(Proposition~\ref{prop:enhanced-sep-limit}).

\subsection{Jacobi forms and Fourier-Jacobi coefficients}

A regular function of index \(m\) on \(\T^\times\) is of the form \(h^m p\), with \(p\) a regular function on \(\T_0\). Its pullback may have poles along the elliptic lattice or at the cusp. For a polynomial in \(a,b,c,t_1,t_2,t_3,h\) homogeneous in weight and index, the pullback is a quasi Jacobi form in the sense of \cite[Appendix~B]{Oberdieck} precisely when it is holomorphic at \(z=0\) for generic \(\tau\) (Theorem~\ref{prop:coordinate-algebra}). The following description of Fourier--Jacobi coefficients follows from the structure theorem for weak Jacobi forms \cite[Theorems~9.3 and~9.4]{EZ} and the classical Fourier--Jacobi theorem \cite{EZ,vdG}. It applies independently of the smoothing family above.

\begin{corollary}\label{thm:intro-FJ}
Let \(\mathfrak F\) be a scalar holomorphic Siegel modular form of integral weight \(r\) for \(\Sp_4(\Z)\) with
trivial character. Write \(q_2=e^{2\pi i\omega}\).
\begin{enumerate}
\setlength{\itemsep}{0pt}
\item Its Fourier-Jacobi coefficients \(\phi_m\in J_{r,m}\) satisfy
\begin{equation}\label{eq:intro-FJ}
\mathfrak F\begin{pmatrix}\tau&z\\z&\omega\end{pmatrix}
=\sum_{m\ge0}\phi_m(\tau,z)q_2^m
=\sum_{m\ge0}p_m(a,b,t_2,t_3)(h q_2)^m,
\end{equation}
where the second equality uses the pullbacks along \(\mathsf t\) on \(\mathcal H\).
\item For weights \(\wt(a,b,t_2,t_3)=(2,3,4,6)\), \(p_m\) has weight \(r+2m\). Modulo \(b^2=4a^3-t_2a-t_3\), it
has the unique expression
\[
p_m=\sum_{j=0}^{m} A_{m,j}a^j+b\sum_{j=0}^{m-2}B_{m,j}a^j,
\qquad A_{m,j},B_{m,j}\in\C[t_2,t_3],
\]
with empty sums interpreted as zero.
\end{enumerate}
\end{corollary}
By Proposition~\ref{prop:weak-ring}, \(\phi_m=\mathsf t^*(h^m p_m)\). Thus the exponent of \(h\) in the regular function \(h^mp_m\) is the Jacobi index of its pullback.
Proposition~\ref{prop:classical-leading} computes the first coefficients of the Igusa cusp forms and the logarithmic
derivatives of the coefficient of \(\widehat\chi_{10}\).

\subsection{Organization}

Sections~\ref{sec:cohom} and~\ref{sec:moduli} describe the limiting cohomology, its cup product, and the four moduli
spaces. Sections~\ref{sec:connection} and~\ref{sec:vector-fields} compute the connection and its vector fields.
Section~\ref{sec:periods} proves that the period map is a biholomorphism; Section~\ref{sec:Jacobi} then determines the
coordinate pullbacks. Section~\ref{sec:degeneration-new} compares the elliptic and genus two period maps along the given smoothing and computes the leading Fourier-Jacobi term. At a separating node, the limiting connection is the direct sum of the two elliptic connections, and its matrix is block diagonal in a split frame.

\section{Relative and limiting de Rham cohomology}\label{sec:cohom}
We consider the curve $C_0=X/(O\sim P)$. Its node is nonseparating,
since its normalization $X$ is connected. We describe the limiting
cohomology, its cup product, and tangentially regularized periods.

\subsection{The relative and open spaces}
Put \(Y=\{O,P\}\) and \(U=X\setminus Y\). The relative de Rham complex is
\(\operatorname{Cone}(\Omega_X^\bullet\to\cO_Y)[-1]\), with \(\cO_Y\) concentrated in degree zero. The logarithmic
complex \(\Omega_X^\bullet(\log Y)\) computes the cohomology of \(U\). The following exact sequences identify the graph
class and the residue quotient.
\begin{lemma}\label{prop:relative-open}
There are exact sequences
\begin{align}\label{eq:relative}
0&\longrightarrow\C\longrightarrow H^1_\dR(X,Y)
 \longrightarrow H^1_\dR(X)\longrightarrow0,\\\label{eq:open}
0&\longrightarrow H^1_\dR(X)\longrightarrow H^1_\dR(U)
 \xrightarrow{\Res}\{(r_O,r_P):r_O+r_P=0\}\longrightarrow0.
\end{align}
Both middle terms have dimension three.
\end{lemma}
\begin{proof}
The long exact sequence for \((X,Y)\) gives \eqref{eq:relative}: the difference of the endpoint values identifies
\(\C^2/\operatorname{diag}\C\) with \(\C\). Accordingly, \(\dd f\) has relative period \(f(P)-f(O)\) if \(f\) is regular
at the endpoints and along the path. For \eqref{eq:open}, the logarithmic residue sequence identifies the quotient of
\(H^1_\dR(U)\) by \(H^1_\dR(X)\) with the residue pairs of sum zero. A differential of the third kind with residues
\((-1,1)\) realizes a generator.
\end{proof}

\subsection{The logarithmic de Rham complex at a node}
The logarithmic de Rham complex identifies the graph, elliptic, and
residue parts of the four dimensional limiting cohomology.

Let \(\nu:X\to C_0=X/(O\sim P)\) be the normalization. Choose a semistable smoothing \(\pi:\mathcal C\to\mathbb D\) and
local coordinates \(\rho_O,\rho_P\) on the two branches at the node, with base coordinate \(s\), in which the smoothing
has equation
\begin{equation}\label{eq:plumb}
 \rho_O\rho_P=s.
\end{equation}
The tangent \(\partial_s|_0\) corresponds to \(\partial_{\rho_O}|_O\otimes\partial_{\rho_P}|_P\in T_OX\otimes T_PX\).
Restriction of the relative logarithmic de Rham complex to the central fiber gives
\begin{equation}\label{eq:log-complex}
 K^\bullet=\left.
 \frac{\Omega^\bullet_{\mathcal C}(\log C_0)}
 {\pi^*\Omega_{\mathbb D}^1(\log0)\wedge\Omega^{\bullet-1}_{\mathcal C}(\log C_0)}
 \right|_{C_0}
 \simeq[\cO_{C_0}\xrightarrow{\dd}\omega_{C_0}].
\end{equation}
Here $\omega_{C_0}$ is the dualizing sheaf. At the node, whose analytic
local ring is $\C\{r,t\}/(rt)$, it is
generated by $e=\dd r/r=-\dd t/t$, and
$\dd f=(r\partial_r f-t\partial_t f)e$. If $q$ denotes the node, the normalization descriptions are
\[
 \begin{aligned}
 \cO_{C_0}
   =\ker\bigl(\nu_*\cO_X
       \xrightarrow{\,f\mapsto f(O)-f(P)\,}\C_q\bigr), 
 \omega_{C_0}
  =\ker\bigl(\nu_*\Omega_X^1(\log Y)
       \xrightarrow{\,\Res_O+\Res_P\,}\C_q\bigr).
 \end{aligned}
\]
Set \(H^1_{\lim,\dR}=\mathbb H^1(C_0,K^\bullet)\), with Hodge filtration induced by \(F^1K^\bullet=[0\to\omega_{C_0}]\).
To identify this complex vector space with the cohomology of a nearby fiber and recover its rational structure and
monodromy, one uses the logarithmic comparison theorem for the chosen smoothing \cite{Steenbrink,Hain}.
\begin{proposition}\label{prop:LMHS}
The limiting de Rham space has dimension four, with \(\dim F^1=2\).
\begin{enumerate}
\setlength{\itemsep}{0pt}
\item Its weight filtration satisfies
\begin{equation}\label{eq:weights}
\begin{gathered}
W_0\simeq\C(0),\qquad W_1/W_0\simeq H^1_\dR(X),\qquad
W_2/W_1\simeq\C(-1),\\
W_1\simeq H^1_\dR(X,Y),\qquad W_2/W_0\simeq H^1_\dR(U).
\end{gathered}
\end{equation}
\item Under logarithmic comparison, \(N=\log T\) has rank one,
\(N^2=0\), \(\im N=W_0\), and \(\ker N=W_1\).
\end{enumerate}
\end{proposition}
\begin{proof}
The normalization sequence and Serre duality give
$h^1(\cO_{C_0})=h^0(\omega_{C_0})=2$.
The logarithmic Hodge-to-de Rham spectral sequence
\(E_1^{p,q}=H^q(C_0,K^p)\Rightarrow\mathbb H^{p+q}(C_0,K^\bullet)\)
degenerates at \(E_1\) \cite{Steenbrink}. Hence \(\dim F^1=2\) and
\(\dim H^1_{\lim,\dR}=2+2=4\). Put
$K^\bullet_{\mathrm{inv}}=[\cO_{C_0}\to\nu_*\Omega_X^1]$.
The residue sequence is
\begin{equation}\label{eq:inv-complex}
0\longrightarrow K^\bullet_{\mathrm{inv}}
 \longrightarrow K^\bullet
 \xrightarrow{\Res_P}\C_q[-1]\longrightarrow0.
\end{equation}
Away from the node, the analytic Poincar\'e lemma applies to \(K^\bullet_{\mathrm{inv}}\). At the node, a holomorphic
differential on each branch has a primitive vanishing at the preimage of the node. These primitives agree there and
hence define a germ in \(\cO_{C_0,q}\); the kernel of \(\dd\) consists of constants. Thus \(K^\bullet_{\mathrm{inv}}\)
resolves \(\C_{C_0}\) analytically. Choose a differential of the third kind on \(X\) with simple poles only at \(O,P\)
and residues \((-1,1)\), as in Lemma~\ref{prop:relative-open}. Its residues sum to zero, so it defines a section of
\(\omega_{C_0}\) and its class maps to \(1\) under \(\Res_P\). By \eqref{eq:inv-complex}, the residue kernel is
\(H^1(C_0,\C)\simeq H^1(X,Y;\C)\).

The classical weight description for a one loop dual graph
\cite[Section~5.4]{Hain} gives
\begin{equation}\label{eq:weight-graded}
 \Gr_0^W H^1_{\lim}\simeq\C(0),\qquad
 \Gr_1^W H^1_{\lim}\simeq H^1(X),\qquad
 \Gr_2^W H^1_{\lim}\simeq\C(-1).
\end{equation}
The invariant cycle theorem identifies the residue kernel with \(W_1\); the graph class spans \(W_0\). Under
restriction to \(U\), this class vanishes, while \(W_1/W_0\) maps isomorphically to \(H^1_\dR(X)\). The chosen
differential maps to \((-1,1)\), a generator of the residue quotient in \eqref{eq:open}. Thus restriction is
surjective, with kernel \(W_0\). The Picard-Lefschetz formula at a nonseparating node gives \(N^2=0\) and
\(\operatorname{rank}N=1\). The monodromy weight filtration then gives \(\im N=W_0\) and \(\ker N=W_1\).
\end{proof}

\subsection{The structure of the cup product}



We next determine the polarization on limiting cohomology. Since the
relevant classes are represented in the \v{C}ech-de Rham total complex,
their transition functions make an essential contribution to the cup
product.

On an acyclic cover, we represent a degree one hypercohomology class by
\((\{\omega_i\},\{f_{ij}\})\), where
\(\omega_j-\omega_i=\dd f_{ij}\) and \(f_{ij}+f_{jk}=f_{ik}\).
Adding \((\{\dd g_i\},\{g_j-g_i\})\) represents the same class.
For the rational expressions in \eqref{eq:frame}, local primitives
remove the auxiliary principal parts, and the resulting data are recorded
in the transition functions.

The complex \(K^\bullet\) carries the usual product of functions with
logarithmic one forms. Although the product of two one forms on a curve
vanishes, the product of two degree one \v{C}ech-de Rham classes need
not vanish: their transition functions contribute in degree two. Using
the cocycle convention above, write
\(\alpha=(\{\omega_i\},\{f_{ij}\})\) and
\(\beta=(\{\nu_i\},\{g_{ij}\})\).
\begin{proposition}\label{prop:cech-cup}
The cup product of these two classes is represented by
\begin{equation}\label{eq:cech-cup}
 (\alpha\cup\beta)^{1,1}_{ij}
   =f_{ij}\nu_j-\omega_i g_{ij},\qquad
 (\alpha\cup\beta)^{2,0}_{ijk}
   =f_{ij}g_{jk}.
\end{equation}
The resulting hypercohomology class is independent of the chosen representatives.
\end{proposition}
\begin{proof}
Multiplication in the double complex gives the displayed components, with the sign $(-1)^{qr}$ when a term of form degree $q$ crosses one of \v Cech degree $r$. Changing either cocycle by a total coboundary changes the product by a total coboundary.
\end{proof}

Put $\chi_{ij}=f_{ij}\nu_j-\omega_i g_{ij}$. To compute the trace, take \(U_0\) to be the complement of a finite set containing the node, and take the other members of the cover to be pairwise disjoint neighborhoods of the omitted points. There are no nonempty triple intersections on this cover.  For \(Q\in X\) with \(\nu(Q)\notin U_0\), let \(i(Q)\) be the unique index \(i>0\) such that \(\nu(Q)\in U_i\). Then
\begin{equation}\label{eq:cup-trace}
 \langle\alpha,\beta\rangle
   =-\sum_Q\Res_Q\chi_{0,i(Q)}.
\end{equation}
At the node, the sum includes the residues at both $O$ and $P$ on the
normalization. This is the Serre duality trace of the complete cocycle
\eqref{eq:cech-cup}, computed on the specified cover.

\begin{example}\label{ex:limiting-cup-matrix}
On \(a\ne0\), set \(f=(x-a)/x\); thus \(f(O)=1\), \(f(P)=0\). The graph class \(\alpha_0=\dd f\) has zero local
differentials and transition constants \(f(O),f(P)\) on the two branches. Since \(\eta_P=(y+b)\dd x/(2(x-a)y)\)
has residues \((-1,1)\), the trace gives
\begin{equation}\label{eq:graph-residue}
\langle\alpha_0,\eta_P\rangle
=-f(O)\Res_O\eta_P-f(P)\Res_P\eta_P=1.
\end{equation}
Thus \(\dd f\) is nonzero in relative cohomology.
Let \(u\) be the primitive of \(\alpha_1\) with \(u(O)=0\). Removing the principal primitive \(-u^{-1}\) of \(x\dd
x/y\) with zero finite constant gives \(\beta\) in \eqref{eq:companion}, with transition \(u^{-1}\) at \(O\).
Hence \(\langle\alpha_1,\beta\rangle=1\). The expansion
\[
\eta_P=\bigl(-u^{-1}-au+\tfrac b2u^2+O(u^3)\bigr)\dd u
\]
gives \(\Res_O(u^{-1}\eta_P)=0\), so \(\langle\beta,\eta_P\rangle=0\).
In contrast, \(\beta_{\mathrm{raw}}=x\dd x/y-\dd(y/(2x))
=\beta+(b/(2a))\alpha_0\) has pairing \(b/(2a)\) with \(\eta_P\).
The classes \(\alpha_1\) and \(\eta_P\) are represented by global sections of \(\omega_{C_0}\), so their cup product
is zero. The graph class \(\alpha_0\) pairs trivially with \(\alpha_1\) and \(\beta\), since both have zero residues.
By skew symmetry the pairing matrix in the ordered basis \((\alpha_0,\alpha_1,\beta,\eta_P)\) is
\[
\Phi=\begin{pmatrix}0&0&0&1\\0&0&1&0\\0&-1&0&0\\-1&0&0&0\end{pmatrix}.
\]
\end{example}

\subsection{Tangential regularization}
An integral of a differential of the third kind diverges logarithmically when the path ends at a pole. We define its
finite part for fixed tangents at \(O\) and \(P\) and determine how it changes when the tangents are rescaled.


Let \(\gamma:[0,1]\to X\) be a piecewise smooth path from \(O\) to \(P\) with interior in \(U\). For \(Q=O,P\) choose a
nonzero tangent \(v_Q\), a local parameter \(\rho_Q\) with \(\dd\rho_Q(v_Q)=1\), and a branch of \(\log\rho_Q\) along
the corresponding end of \(\gamma\). For \(s>0\) close to \(0\) and \(t<1\) close to \(1\), set
\(O_{\eps_O}=\gamma(s)\), \(P_{\eps_P}=\gamma(t)\), where \(\eps_O=\rho_O(\gamma(s))\) and
\(\eps_P=\rho_P(\gamma(t))\). The integral below is taken along \(\gamma|_{[s,t]}\). If \(\omega\) has at most simple
poles at \(O,P\), of residues \(r_O,r_P\), and no poles along the interior of \(\gamma\), define
\begin{equation}\label{eq:reg}
 I_{\rho_O,\rho_P}(\omega)
 =\lim_{\eps_O,\eps_P\to0}
 \left(
   \int_{O_{\eps_O}}^{P_{\eps_P}}\omega
   +r_O\log\eps_O-r_P\log\eps_P
 \right).
\end{equation}
Local primitives \(r_Q\log\rho_Q+g_Q(\rho_Q)\) show that the
two variable limit exists, with \(\eps_O\) and \(\eps_P\) tending to
zero independently, and that the limit is independent of the higher jets
of coordinates normalized by the same tangent. Changing the winding
around \(O\) or \(P\), or changing the branches of the logarithms, can
change the value by integer multiples of \(2\pi i\) times the residues.
We fix these choices whenever comparing periods.

\begin{lemma}\label{lem:tangents}
Let $\rho'_O=A\rho_O+O(\rho_O^2)$ and
$\rho'_P=B\rho_P+O(\rho_P^2)$, where $A,B\in\C^*$.
With compatible logarithm branches,
\[
 I_{\rho'_O,\rho'_P}(\omega)-I_{\rho_O,\rho_P}(\omega)
       =r_O\log A-r_P\log B.
\]
Consequently, for residues $(-1,1)$, replacing $(v_O,v_P)$ by $(\lambda v_O,\mu v_P)$ adds $\log\lambda+\log\mu$.
\end{lemma}
\begin{proof}
Substitute the new parameters in \eqref{eq:reg}. Their logarithmic ratios tend to $\log A,\log B$. Tangent scaling replaces normalized parameters to first order by $\rho_O/\lambda,\rho_P/\mu$.
\end{proof}

For residues $(-1,1)$, inverse tangent rescalings preserve the regularized
value with compatible logarithms. Its class modulo $2\pi i\Z$ therefore
depends on the endpoint tangents only through $v_O\otimes v_P$.

\begin{definition}
\label{def:normalized-tangents}
Let $\alpha_1$ be the distinguished nonzero holomorphic differential.  The
\emph{normalized tangential basepoints} are
$\partial_O=(\alpha_1|_O)^\vee$ and
$\partial_P=(\alpha_1|_P)^\vee$.
A local parameter $\rho_Q$ is normalized when
$\rho_Q(Q)=0$ and $\dd\rho_Q(\partial_Q)=1$.
\end{definition}

In the Weierstrass model these tangents are
\begin{equation}\label{eq:normalized-tangents}
 \partial_O=\partial_{-2x/y},\qquad
 \partial_P=
 \begin{cases}
  b\,\partial_{x-a},&b\ne0,\\[2mm]
  \dfrac{12a^2-t_2}{2}\,\partial_y,&b=0.
 \end{cases}
\end{equation}
The second formula shows that the normalization extends through the branch
points of the Weierstrass projection.


Fix the marked curve \((X,O,P)\), its relative frame \((\alpha_0,\alpha_1,\alpha_2)\), a homology marking, the path
from \(O\) to \(P\), and compatible logarithm branches. In a frame satisfying
Definition~\ref{def:full-T}, write
\(\alpha_3=\eta_P+c\alpha_1+d_T\alpha_0\) and take
\(v_O=\lambda\partial_O,\ v_P=\mu\partial_P\) with \(\lambda,\mu\in\C^\times\).
Let \(K\) be the regularized integral of \(\eta_P+c\alpha_1\) from \(O\) to \(P\) with normalized endpoint
tangents.
\begin{proposition}\label{prop:cup-tangent-calibration}
The closed periods are independent of \(\lambda,\mu\). The regularized relative period of \(\alpha_3\) is
\(K-\mathfrak d\), where
\begin{equation}\label{eq:calibrated-additive-coordinate}
\mathfrak d=d_T-\log\lambda-\log\mu.
\end{equation}
Modulo \(2\pi i\Z\), the relative period is \(K+\log h\), where
\begin{equation}\label{eq:calibrated-multiplicative-coordinate}
h=e^{-\mathfrak d}=\lambda\mu e^{-d_T},
\end{equation}
which is independent of the logarithm branches and invariant under
\((\lambda,\mu)\mapsto(r\lambda,r^{-1}\mu)\) for \(r\in\C^\times\).
\end{proposition}

\begin{proof}
Closed periods involve only cohomology classes. By Lemma~\ref{lem:tangents}, rescaling the tangents adds
\(\log\lambda+\log\mu\) to \(K\), since the residues are \((-1,1)\). The term \(d_T\alpha_0\) contributes \(-d_T\)
because the relative period of \(\alpha_0\) is \(-1\). This yields \(K-\mathfrak d\); modulo \(2\pi i\Z\) it
equals \(K+\log h\). The two invariances follow from \eqref{eq:calibrated-multiplicative-coordinate}.
\end{proof}

\section{Enhanced elliptic curves and the additional coordinate}\label{sec:moduli}

\subsection{The enhanced elliptic curve with two marked points}

The space \(\T_0\) of \cite{CMVL} parametrizes elliptic curves with two marked points and a relative de Rham frame. We
describe this frame in Weierstrass coordinates and construct its class of the second kind also when \(a=0\).

Recall that an enhanced elliptic curve with two marked points consists
of a triple \((X,O,P)\) with \(P\ne O\), a nonzero holomorphic
differential \(\alpha_1\), and a frame
\((\alpha_0,\alpha_1,\alpha_2)\) of \(H^1_{\mathrm{dR}}(X,Y)\) such that
\(\alpha_0\) is the graph class whose relative period is \(-1\), and the images of \(\alpha_1,\alpha_2\) in
\(H^1_{\dR}(X)\) have cup product \(1\). The elliptic curve carries its principal polarization. This is
the index zero moduli problem represented by \(\T_0\) in \cite{CMVL}.
The relative class \(\alpha_0\) is canonical once the ordered pair
\((O,P)\) and the period normalization are fixed.


The origin \(O\) and differential \(\alpha_1\) determine the Weierstrass coordinates of \eqref{eq:curve}, in which
\(\alpha_1=\dd x/y\).
Writing \(P=(a,b)\) gives \(t_3=4a^3-t_2a-b^2\).
Let \(\beta\) be the relative class of \(x\dd x/y\), obtained by subtracting at \(O\)
the differential of a local primitive of its polar part, with constant term zero.
In absolute cohomology,
\(\langle\alpha_1,\beta\rangle=\langle\alpha_1,\alpha_2\rangle=1\).
Thus \(\alpha_2-\beta\) lies in the span of \(\alpha_0,\alpha_1\). Thus there are unique \(c,t_1\) such that
\begin{equation}\label{eq:companion}
 \alpha_2=c\alpha_0+t_1\alpha_1+\beta.
\end{equation}
On $a\ne0$,
\[
 \beta=\frac{x\dd x}{y}-\dd\left(\frac y{2x}\right)-\frac b{2a}\alpha_0.
\]
This rational expression is used only on \(a\ne0\). For every \(a\), the
function \(\widetilde f=(x-a)/(x-a+1)\) has endpoint values \((1,0)\).
Its differential, together with the local primitives at its auxiliary poles,
represents \(\alpha_0\). Let \(u\) satisfy \(\dd u=\alpha_1\) and \(u(O)=0\).
Subtracting the differential of the principal primitive \(-u^{-1}\) at \(O\),
with zero finite constant and no subtraction at \(P\), defines \(\beta\)
also at \(a=0\).

\subsection{The residue direction and polarization}

We choose a differential of the third kind with residues \(-1\) at \(O\) and \(1\) at \(P\), and compute its pairings with the relative frame. These pairings determine the symplectic extensions of that frame.

Set
\[
  \eta_P=\frac{y+b}{2(x-a)}\frac{\dd x}{y}.
\]
It has simple poles at \(O,P\) with respective residues \(-1,1\), and hence defines a section of \(\omega_{C_0}\). If
\(b\ne0\), its apparent singularity at
\((a,-b)\) is removable. If \(b=0\), this point is \(P\);
using \(y\) as a local parameter at \(P\), one finds that
\(\eta_P-\dd y/y\) is holomorphic there.
The tangent dual to \(\alpha_1\) remains nonzero at \(b=0\), as in
\eqref{eq:normalized-tangents}.

\begin{proposition}\label{prop:polarized-frame}
For \((e_0,e_1,e_2,e_3)=(\alpha_0,\alpha_1,\beta,\eta_P)\), the limiting pairing satisfies
\begin{equation}\label{eq:Phi}
\bigl(\langle e_i,e_j\rangle\bigr)_{0\le i,j\le3}=\Phi,
\end{equation}
with \(\Phi\) as in Example~\ref{ex:limiting-cup-matrix}.
\end{proposition}
\begin{proof}
Example~\ref{ex:limiting-cup-matrix} proves the identity on \(a\ne0\) using \eqref{eq:cup-trace}. At \(a=0\), use
\(\dd\widetilde f\) for the graph class. Since
\(\widetilde f(O)=1\) and \(\widetilde f(P)=0\), the same residue
computation applies.
\end{proof}
A symplectic extension of \eqref{eq:companion}, with zero pairing against $\alpha_1,\alpha_2$, is therefore
$\alpha_3=\eta_P+c\alpha_1+d\alpha_0$.
The matrix sending the base frame to this frame is
\[
 S=
 \begin{pmatrix}
  1&0&0&0\\
  0&1&0&0\\
  c&t_1&1&0\\
  d&c&0&1
 \end{pmatrix},
 \qquad S\Phi S^{\tr}=\Phi.
\]
An arbitrary such $\alpha_3$ need not belong to $F^1$. The de Rham
Hodge filtration is induced by $K^\bullet$; its transport to fixed
nearby cohomology depends on the smoothing comparison. The weight and
pairing conditions alone do not specify $F^1$.

On the open set \(a\ne0\), the frame is given by
\begin{align}\label{eq:frame}
 \alpha_0=&\dd\left(\frac{x-a}{x}\right),
& \alpha_1&=\frac{\dd x}{y},\notag\\
 \alpha_2=&\left(c-\frac b{2a}\right)\alpha_0+t_1\frac{\dd x}{y}+\frac{x\dd x}{y}-\dd\left(\frac y{2x}\right),
 &\alpha_3&=\frac{y+b}{2(x-a)}\frac{\dd x}{y}+c\alpha_1+d\alpha_0.
\end{align}

These expressions represent cohomology classes whose \v{C}ech data
include local primitives at the auxiliary poles. In particular,
\(\alpha_0\) is the nonzero graph class. The intrinsic construction of
\(\beta\) above extends the frame through \(a=0\).

\subsection{The full T-space of tangentially enhanced frames}\label{subsec:full-T}

A point of \(\T\) includes a choice of the class \(\alpha_3\) and a nonzero tangent vector at each marked point. Fixing
the two tangents dual to \(\alpha_1\) gives the section \(\T^{\log}\). In the next subsection we use their tensor
product to construct \(\T^\times\).

\begin{definition}\label{def:full-T}
A \emph{tangentially enhanced elliptic frame} consists of a principally
polarized elliptic curve \((X,O)\), a point \(P\ne O\), nonzero tangents
\(v_O\in T_OX\), \(v_P\in T_PX\), and an ordered basis
\((\alpha_0,\alpha_1,\alpha_2,\alpha_3)\) of \(H^1_{\lim,\dR}\) for
\(X/(O\sim P)\), subject to the following conditions.
\begin{enumerate}
\item \(\alpha_0\) is the graph class in \(W_0\) with period \(-1\) along a relative path from \(O\) to \(P\);
\(\alpha_1\in F^1\cap W_1\) and \(\alpha_2\in W_1\).
\item \(\alpha_3\) has residues \(-1\) at \(O\) and \(1\) at \(P\), and
\((\langle\alpha_i,\alpha_j\rangle)_{i,j}=\Phi\).
\end{enumerate}
We write \(\T\) for the space of isomorphism classes. An isomorphism carries \(O\) to \(O'\), \(P\) to \(P'\), and
\(v_Q\) to \(v'_{Q'}\) under its differential at each marked point; it pulls back each target frame class to the
corresponding source class.
\end{definition}
The images of \(\alpha_1,\alpha_2\) then form a symplectic basis of
\(H^1_{\dR}(X)\); no condition \(\alpha_3\in F^1\) is imposed.
Here \(H^1_{\lim,\dR}=\mathbb H^1(C_0,K^\bullet)\), with \(K^\bullet\) as in
\eqref{eq:log-complex}. Its weight filtration and pairing are those of Section~\ref{sec:cohom}.
A choice of smoothing and a comparison with nearby cohomology enter only when we consider
the limiting rational mixed Hodge structure and its periods.

\begin{theorem}\label{thm:T-representation}
Tangentially enhanced elliptic frames over \(\C\) have unique normalized
Weierstrass representatives, parametrized by the smooth eight dimensional variety
\[
\T\simeq\Spec A_0[d_T,\lambda^{\pm1},\mu^{\pm1}],
\qquad A_0=\C[a,b,c,t_1,t_2,\Delta^{-1}].
\]
The curve is \eqref{eq:curve}; the frame is \eqref{eq:frame} with \(d=d_T\),
using its intrinsic extension at \(a=0\). The tangents are
\begin{equation}\label{eq:full-tangents}
v_O=\lambda(\alpha_1|_O)^\vee,\qquad v_P=\mu(\alpha_1|_P)^\vee.
\end{equation}
Forgetting the tangents gives \(\T\to\T^{\log}=\Spec A_0[d]\),
with \(d\mapsto d_T\) on coordinate rings and section \(\lambda=\mu=1\).
\end{theorem}
\begin{proof}
The differential $\alpha_1$ determines the normalized Weierstrass form
$y^2=4x^3-t_2x-t_3$, $\alpha_1=\dd x/y$. The point $P=(a,b)$ gives
$t_3=4a^3-t_2a-b^2$, and smoothness is equivalent to $\Delta\ne0$.
The relative frame expression \eqref{eq:companion} determines $c,t_1$.

Write $\alpha_3=\eta_P+u_0\alpha_0+u_1\alpha_1+u_2\beta$.
The pairing relations \eqref{eq:Phi} give $u_2=0$ and $u_1=c$;
the coefficient $u_0=d_T$ is unrestricted. Conversely,
$S\Phi S^{\tr}=\Phi$ verifies all the required pairings.

Finally, $\alpha_1$ trivializes each tangent line and determines
$\lambda,\mu$ uniquely. At $b=0$,
$\dd x/y=2\dd y/(12x^2-t_2)$ and $12a^2-t_2\ne0$ by smoothness,
so these trivializations extend across the branch points.
Forgetting $\lambda,\mu$ gives the stated projection, and setting
$\lambda=\mu=1$ gives its canonical section.
\end{proof}

\subsection{The effective multiplicative enhancement}

By Proposition~\ref{prop:cup-tangent-calibration}, the regularized period depends
on \(d_T,\lambda,\mu\) only through \(\lambda\mu e^{-d_T}\) when the period is taken modulo \(2\pi i\Z\).
We now identify the tangent tensor represented by this scalar.

\begin{definition}\label{def:effective-T}
For a point of $\T$, define its \emph{effective tangent tensor} by
\begin{equation}\label{eq:effective-tensor}
 \upsilon_{\mathrm{eff}}
 :=k^2e^{-d_T}v_O\otimes v_P\in T_OX\otimes T_PX.
\end{equation}
Two tangentially enhanced frames have the same effective enhancement when
their relative data $(X,O,P;\alpha_0,\alpha_1,\alpha_2)$ agree and their
effective tangent tensors agree.  The parameter space of effective
enhancements is denoted by $\T^\times$.
\end{definition}


\begin{remark}
Over a fixed relative enhanced elliptic curve, the effective enhancements
are the nonzero tensors in \(T_OX\otimes T_PX\). Indeed, by
\eqref{eq:full-tangents},
\begin{equation}\label{eq:tensor}
 \upsilon_{\mathrm{eff}}=k^2h\,(\alpha_1|_O)^\vee\otimes
 (\alpha_1|_P)^\vee,\qquad h=\lambda\mu e^{-d_T}.
\end{equation}
Every \(h\ne0\) occurs: take \((d_T,\lambda,\mu)=(0,h,1)\).
Thus the tensor description gives \(\T^\times=\Spec A_0[h,h^{-1}]\).
The map from \(\T^{\mathrm{an}}\) to \(\T^{\times,\mathrm{an}}\) is
holomorphic but not algebraic: on the section \(\lambda=\mu=1\),
\begin{equation}\label{eq:effective}
 d=d_T,\qquad h=e^{-d}.
\end{equation}
Replacing a local logarithm \(d=-\log h\) by \(d+2\pi i n\) adds
\(2\pi i n\alpha_0\) to \(\alpha_3\) and leaves \(h\) unchanged.
For the descended connection, see Corollary~\ref{prop:geometric}.
\end{remark}

\subsection{Changes of frame}
We compute how a change of de Rham frame acts on the full space \(\T\) and on the effective tensor.

Let \(\G\) be the set of matrices
\begin{equation}\label{eq:G}
g(s,s',v,w)=
\begin{pmatrix}
1&0&v&w\\ 0&s&s'&sv\\ 0&0&s^{-1}&0\\ 0&0&0&1
\end{pmatrix},
\qquad s\in\C^\times,\quad s',v,w\in\C.
\end{equation}
\begin{proposition}\label{prop:frame-change}
The stabilizer of \(\alpha_0\), the weight filtration,
\(F^1\cap W_1=\C\alpha_1\), and \(\Phi\) is \(\G\).
Right multiplication of the row \((\alpha_0,\alpha_1,\alpha_2,\alpha_3)\) by \(g(s,s',v,w)\)
induces the following algebraic action on \(\T\):
\[
\begin{alignedat}{4}
a'&=s^{-2}a, &\quad b'&=s^{-3}b, &\quad c'&=v+s^{-1}c,
&\quad d_T'&=d_T+w,\\
t_1'&=s^{-1}s'+s^{-2}t_1, & t_2'&=s^{-4}t_2,
& \lambda'&=s\lambda, & \mu'&=s\mu.
\end{alignedat}
\]
On effective enhancements,
\begin{equation}\label{eq:effective-action}
h'=s^2e^{-w}h,\qquad d'=d+w-2\log s
\quad\text{on compatible logarithmic charts for \(d\).}
\end{equation}
\end{proposition}
\begin{proof}
The flag and the pairing on \(W_1/W_0\) give
\(\alpha'_1=s\alpha_1\) and
\(\alpha'_2=v\alpha_0+s'\alpha_1+s^{-1}\alpha_2\).
Write \(\alpha'_3=A\alpha_0+B\alpha_1+C\alpha_2+D\alpha_3\).
Its pairings with \(\alpha_0,\alpha'_1,\alpha'_2\) give, respectively,
\(D=1\), \(C=0\), and \(B=sv\). Taking \(w=A\) gives \eqref{eq:G}. Its action is
\[
 (\alpha_0,\alpha_1,\alpha_2,\alpha_3)g
 =\bigl(\alpha_0,\ s\alpha_1,
        v\alpha_0+s'\alpha_1+s^{-1}\alpha_2,
        w\alpha_0+sv\alpha_1+\alpha_3\bigr).
\]
Set $x'=s^{-2}x$, $y'=s^{-3}y$. The Weierstrass equation and
\eqref{eq:companion} give the displayed coordinate transformations;
the normalized second kind class scales by $s^{-1}$ and $\eta_P$
is invariant. Keeping the endpoint tangents fixed gives
$\lambda'=s\lambda$, $\mu'=s\mu$, since
$(\alpha'_1|_Q)^\vee=s^{-1}(\alpha_1|_Q)^\vee$.
Substitution in $h=\lambda\mu e^{-d_T}$ and a local logarithm prove
\eqref{eq:effective-action}.
\end{proof}

\begin{remark}
If the endpoint tangents are held fixed, a change of frame sends
\((\lambda,\mu)=(1,1)\) to \((s,s)\). Thus the normalized section is not
preserved in general. The translation \(d_T'=d_T+w\) is algebraic on \(\T\),
whereas the formula for \(d'\) in \eqref{eq:effective-action} requires a local
choice of \(\log s\). On \(\T^\times\) the central action \(h\mapsto e^{-w}h\)
is globally holomorphic but not algebraic: every regular map \(\mathbb A^1\to\Gm\) is constant.
\end{remark}

Let \(\G_0\) be the subgroup of frame changes that fix \(d_T\); by
Proposition~\ref{prop:frame-change}, these are exactly the matrices \(g(s,s',v,0)\).
\begin{proposition}\label{prop:effective-algebraic-group}
The central translations of \(d_T\) give a direct product decomposition
\(\G=\G_0\times\Ga\). The group \(\G_{\mathrm{eff}}=\G_0\times\Gm\) acts algebraically on \(\T^\times\):
\(\G_0\) acts on \(a,b,c,t_1,t_2\) as in
Proposition~\ref{prop:frame-change}, and \(h'=s^2\xi h\) for
\(\xi\in\Gm\). The action of \(\G^{\rm an}\) on \((\T^\times)^{\rm an}\)
is the pullback of this action by the surjective homomorphism
\[
\G^{\rm an}\longrightarrow\G_{\mathrm{eff}}^{\rm an},\qquad
(g_0,w)\longmapsto(g_0,e^{-w}),
\]
whose kernel is the central subgroup \(\{1\}\times 2\pi i\Z\).
\end{proposition}
\begin{proof}
In \eqref{eq:G} the coordinate \(w\) is additive under multiplication,
and the translations \(g(1,0,0,w)\) are central. This gives the direct product.
The transformations of the first five coordinates are those of \(\G_0\);
the factors \(s^2\xi\) multiplying \(h\) are multiplicative.
The coordinate transformations are polynomial in \(s,s^{-1},s',v,\xi,\xi^{-1}\)
and the coordinates of \(\T^\times\). Since
\(\Delta'=s^{-12}\Delta\) and \((h')^{-1}=s^{-2}\xi^{-1}h^{-1}\),
they also define pullbacks of the inverted generators. The kernel and surjectivity
follow from those of \(w\mapsto e^{-w}\).
\end{proof}

\begin{remark}
Differentiating the frame action in the \(s,s',v,w\) directions at the identity
gives \(\R_1,\R_2,\R_3,\R_4\), respectively, with \(\R_4=-h\partial_h\)
on \(\T^\times\). Indeed, \(w\mapsto e^{-w}\) sends \(\partial_w\) to
\(-\xi\partial_\xi\). The two remaining fields, \(\R_\tau\) and \(\R_z\),
are defined by the connection in Section~\ref{sec:vector-fields}.
\end{remark}

\section{The Gauss-Manin connection}\label{sec:connection}
We compute the Gauss-Manin connection in the chosen frame by
differentiating rational representatives at fixed $x$ and reducing in
relative and logarithmic de Rham cohomology.

\subsection{Algebraic reduction and endpoint terms}
We give the reduction identities used to differentiate periods.
In these identities, differentials are taken along the curve fibers,
and a prime denotes differentiation with respect to $x$.
Exact terms are kept because their endpoint values and finite
parts contribute to relative and third kind periods.

\begin{lemma}[{\cite[proof of Proposition~2.2]{Mov}}]
\label{lem:polynomial-reduction}
Let $p(x)=4x^3-t_2x-t_3$. For $n\in\Z_{\ge0}$,
\begin{equation}\label{eq:reduction1}
\dd(x^ny)=\left((4n+6)x^{n+2}
  -\left(n+\tfrac12\right)t_2x^n-t_3(x^n)'\right)\frac{\dd x}{y}.
\end{equation}
Thus every $Q(x)\dd x/y$, with $Q\in\C[x]$, is a linear combination
of $\dd x/y$ and $x\dd x/y$ modulo an exact differential.
\end{lemma}

\begin{example}
Applying \eqref{eq:reduction1} successively for \(n=0,1,2,3\) gives
\begin{equation}\label{eq:high-reductions}
\begin{aligned}
\frac{x^2\dd x}{y}
 &=\frac{t_2}{12}\frac{\dd x}{y}+\dd\!\left(\frac{y}{6}\right),\\
\frac{x^3\dd x}{y}
 &=\frac{t_3}{10}\frac{\dd x}{y}
   +\frac{3t_2}{20}\frac{x\dd x}{y}
   +\dd\!\left(\frac{xy}{10}\right),\\
\frac{x^4\dd x}{y}
 &=\frac{5t_2^2}{336}\frac{\dd x}{y}
   +\frac{t_3}{7}\frac{x\dd x}{y}
   +\dd\!\left(\frac{x^2y}{14}+\frac{5t_2y}{168}\right),\\
\frac{x^5\dd x}{y}
 &=\frac{t_2t_3}{30}\frac{\dd x}{y}
   +\frac{7t_2^2}{240}\frac{x\dd x}{y}+\dd\!\left(\frac{x^3y}{18}
              +\frac{7t_2xy}{360}+\frac{t_3y}{36}\right).
\end{aligned}
\end{equation}
\end{example}

The preceding four identities, including their exact terms, are given in
\cite[\S6.2]{CMVL}. The same section gives the B\'ezout identity and
the case \(Q=1\) of \eqref{eq:reduction2}; the general formula follows
by the product rule. It also gives the moving endpoint term \(\dd a/b\)
in the variation of \(\int_O^P\dd x/y\) for \(b\ne0\).

Set \(p=4x^3-t_2x-t_3\) and define
\[
A_1=-36x^4+15t_2x^2-t_2^2,\qquad A_2=-108x^3+27t_2x-27t_3.
\]
\begin{lemma}\label{lem:cubic-denominator-reduction}
For \(\Delta\ne0\) and \(Q\in\C[x]\),
\begin{equation}\label{eq:reduction2}
\frac{Q(x)\dd x}{y^3}
=\frac{QA_2-2(QA_1)'}{\Delta}\frac{\dd x}{y}
+\dd\!\left(\frac{2QA_1}{\Delta y}\right).
\end{equation}
\end{lemma}
\begin{proof}
Differentiate $2QA_1/(\Delta y)$ and use $2y\dd y=p'(x)\dd x$ together
with the B\'{e}zout identity $-p'A_1+pA_2=\Delta$.  The remaining numerator is
then reduced by Lemma~\ref{lem:polynomial-reduction}.
\end{proof}

Put \(\omega_i=x^i\dd x/y\) for \(i\ge0\), with independent base
coordinates \(a,b,t_2\).
\begin{lemma}\label{lem:parameter-derivatives}
At fixed \(x\), for \(\xi\in\{a,b,t_2\}\),
\begin{equation}\label{eq:raw-derivatives}
(\partial_a p,\partial_b p,\partial_{t_2}p)=(t_2-12a^2,\,2b,\,a-x),
\qquad
\partial_\xi\omega_i=-\frac{x^i\partial_\xi p}{2y^3}\,\dd x .
\end{equation}
For \(b\ne0\), the moving endpoint contributes \(1/b\) to
\(\partial_a\int_O^P\omega_0\) at fixed \(b,t_2\).
\end{lemma}
\begin{proof}
Differentiate \(p=4x^3-t_2x-(4a^3-t_2a-b^2)\) and \(y^{-1}=p^{-1/2}\)
on a local sheet of the Weierstrass projection. Since \(x(P)=a\)
and \(y(P)=b\ne0\), the Leibniz rule gives the endpoint term
\((1/y)(P)\,\partial_a x(P)=1/b\).
\end{proof}

On a simply connected marked chart with \(b\ne0\), fix the path classes
and logarithms. Choose the path interiors to avoid the finite
ramification points. For the third kind form \(\eta_P\) of Section~3.2 put
\[
\mathcal I_{\mathrm W}=
\Reg_{(\partial_{-2x/y},\partial_{x-a})}\int_O^P\eta_P .
\]
\begin{lemma}\label{lem:third-kind-variation}
The derivatives of the regularized period in the coordinates \(a,b,t_2\) are
\begin{equation}\label{eq:third-kind-raw-variation}
\begin{aligned}
\partial_a\mathcal I_{\mathrm W}
&=-3b\int_O^P(x+a)\frac{\dd x}{y^3},\\
\partial_b\mathcal I_{\mathrm W}
&=\int_O^P\left(2(x^2+ax+a^2)-\frac{t_2}{2}\right)\frac{\dd x}{y^3},\\
\partial_{t_2}\mathcal I_{\mathrm W}
&=\frac b4\int_O^P\frac{\dd x}{y^3}.
\end{aligned}
\end{equation}
For the normalized tangents of Definition~\ref{def:normalized-tangents},
\begin{equation}\label{eq:normalized-third-kind-correction}
\mathcal I_{\alpha_1}=\mathcal I_{\mathrm W}+\log b,
\qquad
\dd\mathcal I_{\alpha_1}=\dd\mathcal I_{\mathrm W}+\frac{\dd b}{b}.
\end{equation}
\end{lemma}
\begin{proof}
Set $g(x)=(1+b/y)/(2(x-a))$, so that $\eta_P=g(x)\dd x$.
Differentiation at fixed $x$, using $2y\partial_\xi y=\partial_\xi p$,
gives
\begin{align*}
 \partial_a\eta_P
   =-\frac{3b(x+a)}{y^3}\dd x-\dd_x g, \qquad
 \partial_b\eta_P
   =\frac{2(x^2+ax+a^2)-t_2/2}{y^3}\dd x,\qquad
 \partial_{t_2}\eta_P=\frac{b}{4y^3}\dd x.
\end{align*}
Here \(\dd_x\) denotes differentiation along the fiber. At the upper cutoff \(x=a+\epsilon_P\),
differentiation with respect to \(a\) contributes \(g(a+\epsilon_P)\), which cancels
the endpoint value of the primitive \(-g\). For the other two derivatives this cutoff is fixed.
At the lower cutoff, \(\rho_O=-2x/y\), one has \(x=\rho_O^{-2}+O(\rho_O^2)\)
and \(g=O(\rho_O^2)\), so the lower endpoint terms tend to zero.
The residues are constant, hence the logarithmic subtraction terms have zero derivative
at fixed \(\rho_O,\epsilon_P\). The remaining integrands are regular at both endpoints;
passing to the limit gives \eqref{eq:third-kind-raw-variation}.
At \(P\) the tangent dual to \(\alpha_1\) is \(b\partial_{x-a}\), whereas the tangent at \(O\)
is already normalized. Lemma~\ref{lem:tangents} gives the correction \(\log b\) to the integral
and hence the term \(\dd b/b\) in \eqref{eq:normalized-third-kind-correction}.
\end{proof}

\subsection{The connection matrix in the chosen frame}
We express the derivatives of the four frame classes in terms of six one forms on \(\T^{\log}\).
Writing the first three forms as the solution of a linear system makes their regularity
on \(\Delta\ne0\) apparent.

Let
\[
\begin{aligned}
A_a=-2a^2+2at_1+bc+\frac{t_2}{3},\qquad
A_b=\left(6a^2-\frac{t_2}{2}\right)c-3ab+3bt_1,
\qquad A_t&=4t_1t_2-6t_3,
\end{aligned}
\]
and define
\begin{equation}\label{eq:Nmatrix}
M=\begin{pmatrix}
A_a & b & -2a\\
A_b & 6a^2-t_2/2 & -3b\\
A_t & 0 & -4t_2
\end{pmatrix},
\qquad \det M=-\frac23\Delta.
\end{equation}
Set
\begin{equation}\label{eq:coframe}
\begin{aligned}
\begin{pmatrix}\vartheta_\tau\\\vartheta_z\\\vartheta_1\end{pmatrix}
&=-\frac{3}{2\Delta}\operatorname{adj}(M)
  \begin{pmatrix}\dd a\\\dd b\\\dd t_2\end{pmatrix},\\
\vartheta_2
&=\dd t_1-\left(t_1^2-\frac{t_2}{12}\right)\vartheta_\tau
  +2t_1\vartheta_1,\\
\vartheta_3
&=\dd c-\left((a+t_1)c-\frac{b}{2}\right)\vartheta_\tau
  -(a+t_1)\vartheta_z+c\vartheta_1,\\
\vartheta_4
&=\dd d-(c^2-a+2t_1)\vartheta_\tau
  -2c\vartheta_z+2\vartheta_1.
\end{aligned}
\end{equation}
Since $\det M=-2\Delta/3$, all six forms are regular on $\T^{\log}$,
including at $a=0$ and $b=0$.

\begin{proposition}\label{prop:explicit-connection}
In the ordered frame
$\alpha=(\alpha_0,\alpha_1,\alpha_2,\alpha_3)^{\tr}$,
the Gauss-Manin connection is
\begin{equation}\label{eq:connection}
\nabla\alpha=\A\alpha,
\qquad
\A=\begin{pmatrix}
0&0&0&0\\
-\vartheta_z&\vartheta_1&-\vartheta_\tau&0\\
\vartheta_3&\vartheta_2&-\vartheta_1&0\\
\vartheta_4&\vartheta_3&\vartheta_z&0
\end{pmatrix}.
\end{equation}
The matrix is regular on $\Delta\ne0$ and satisfies
$\A\Phi+\Phi\A^{\tr}=0$.
\end{proposition}
\begin{proof}
The graph class $\alpha_0$ is locally constant, so $\nabla\alpha_0=0$.
Differentiate the representatives of $\alpha_1$ and $\alpha_2$ at fixed $x$,
apply \eqref{eq:reduction2}, and then use \eqref{eq:reduction1} to reduce the
polynomial numerator.  The exact primitives contribute their endpoint
differences in relative cohomology.  Substitution of
$\alpha_2=c\alpha_0+t_1\alpha_1+\beta$ gives the second and third rows of
\eqref{eq:connection}.

Write \(\nabla\alpha_3=u_0\alpha_0+u_1\alpha_1+u_2\alpha_2\); the residue is constant,
so there is no \(\alpha_3\) term. Differentiating
\(\langle\alpha_1,\alpha_3\rangle=\langle\alpha_2,\alpha_3\rangle=0\)
and using the preceding two rows gives \(u_2=\vartheta_z\) and \(u_1=\vartheta_3\).
To determine \(u_0\), let \(I_j\) be the regularized period of \(\alpha_j\) along the chosen path
from \(O\) to \(P\). Since \(I_0=-1\) and \(I_3=\mathcal I_{\alpha_1}+cI_1-d\), we have
\[
u_0=\vartheta_3 I_1+\vartheta_z I_2-\dd I_3.
\]
Use the second row for \(\dd I_1\) and \eqref{eq:normalized-third-kind-correction} for \(\dd\mathcal I_{\alpha_1}\).
Reducing the three integrands in \eqref{eq:third-kind-raw-variation} by \eqref{eq:reduction2} then gives
\(u_0=\vartheta_4\) as in \eqref{eq:coframe}.
This proves \eqref{eq:connection} on \(ab\ne0\). The frame and the displayed matrix are regular
on \(\Delta\ne0\), so the identity extends across \(a=0\) and \(b=0\).
The relation \(\A\Phi+\Phi\A^{\tr}=0\) also follows directly from the matrix.
\end{proof}
\subsection{Flatness and conventions}
For our column of basis elements, \(\nabla\alpha=\A\alpha\). By the graded
Leibniz rule, \(\nabla^2\alpha=(\dd\A-\A\wedge\A)\alpha\).
Thus the flatness equation in this convention is
\begin{equation}\label{eq:flatness}
\dd\A=\A\wedge\A.
\end{equation}
Suppose that \(\A(V)=C_V\) and \(\A(W)=C_W\) are constant matrices.
Then \(\dd\A(V,W)=-\A([V,W])\), whereas
\((\A\wedge\A)(V,W)=C_VC_W-C_WC_V\). Hence the curvature vanishes on
\((V,W)\) if and only if
\begin{equation}\label{eq:anti}
\A([V,W])=[C_W,C_V].
\end{equation}
In Section~\ref{sec:vector-fields} we verify this equality on a basis of vector fields.
To pass to \(\T^\times\), substitute \(\dd d=-\dd h/h\) in \eqref{eq:coframe} and \eqref{eq:connection}
and denote the resulting matrix by \(\A^\times\). On
$\mathcal H^\times=\cO_{\T^\times}^{\oplus4}$ with frame
$\epsilon=(\epsilon_0,\epsilon_1,\epsilon_2,\epsilon_3)^{\tr}$,
define $\nabla^\times\epsilon=\A^\times\epsilon$ and
$\langle\epsilon_i,\epsilon_j\rangle=\Phi_{ij}$.
The matrix is regular, flat, and pairing preserving. Its period
realization is proved in Corollary~\ref{prop:geometric}.
We continue to write $\A$ for $\A^\times$ when the base is $\T^\times$;
the corresponding vector fields are obtained by
$\partial_d\mapsto-h\partial_h$. On $\T$, the effective connection
pulls back analytically by replacing $\dd d$ with
$\dd d_T-\dd\lambda/\lambda-\dd\mu/\mu$.

\section{Modular vector fields}\label{sec:vector-fields}

In this section we prove the vector field and flatness assertions of
Theorem~\ref{thm:main}. We determine the six fields from the connection,
compute their Lie brackets, and verify regularity on the discriminant
complement. Let $E_{ij}$ denote the elementary matrix with row and column
indices $0,1,2,3$, and set
\begin{equation}\label{eq:matrices}
\begin{alignedat}{3}
C_\tau&=-E_{12},&\qquad C_z&=-E_{10}+E_{32},
 &\qquad C_1&=E_{11}-E_{22},\\
C_2&=E_{21},& C_3&=E_{20}+E_{31},
 & C_4&=E_{30}.
\end{alignedat}
\end{equation}

Use \(A_a,A_b,A_t\) from the first column of \(M\) in \eqref{eq:Nmatrix}.
\begin{proposition}\label{prop:fields-from-connection}
On \(\T^{\log}\), the equations \(\A(\R_i)=C_i\),
\(i\in\{\tau,z,1,2,3,4\}\), have unique regular solutions:
\begin{equation}\label{eq:fields}
\begin{aligned}
\R_\tau={}&A_a\partial_a+A_b\partial_b+
\bigl((a+t_1)c-b/2\bigr)\partial_c+(c^2-a+2t_1)\partial_d\\
&+(t_1^2-t_2/12)\partial_{t_1}+A_t\partial_{t_2},\\
\R_z={}&b\partial_a+(6a^2-t_2/2)\partial_b
+(a+t_1)\partial_c+2c\partial_d,\\
\R_1={}&-2a\partial_a-3b\partial_b-c\partial_c-2\partial_d
-2t_1\partial_{t_1}-4t_2\partial_{t_2},\\
\R_2={}&\partial_{t_1},\qquad \R_3=\partial_c,\qquad \R_4=\partial_d .
\end{aligned}
\end{equation}
Their values form a basis of the tangent space at every point of \(\T^{\log}\).
Indeed, in the coordinate order \((a,b,c,d,t_1,t_2)\),
\begin{equation}\label{eq:det}
\det(\R_\tau,\R_z,\R_1,\R_2,\R_3,\R_4)=\frac23\Delta .
\end{equation}
\end{proposition}
\begin{proof}
By \eqref{eq:connection}, the defining condition is equivalent to
$\vartheta_j(\R_i)=\delta_{ij}$.  Since $(\dd a,\dd b,\dd t_2)^{\tr}
=M(\vartheta_\tau,\vartheta_z,\vartheta_1)^{\tr}$, the
$(a,b,t_2)$ components of $\R_\tau,\R_z,\R_1$ are the columns
of $M$; they vanish for $\R_2,\R_3,\R_4$. The remaining
components follow from \eqref{eq:coframe}. This gives
exactly \eqref{eq:fields}.  Taking determinants and using
$\det M=-2\Delta/3$ gives \eqref{eq:det}.
\end{proof}

\begin{proposition}\label{prop:brackets}
The nonzero brackets, up to skew symmetry, are
\begin{equation}\label{eq:brackets}
\begin{alignedat}{4}
\bigl[\R_1,\R_\tau\bigr] &= -2\R_\tau, &\quad
\bigl[\R_1,\R_z\bigr]     &= -\R_z,     &\quad
\bigl[\R_1,\R_2\bigr]     &= 2\R_2,     &\quad
\bigl[\R_1,\R_3\bigr]     &= \R_3,\\[2pt]
\bigl[\R_2,\R_\tau\bigr] &= -\R_1,      &\quad
\bigl[\R_2,\R_z\bigr]     &= \R_3,      &\quad
\bigl[\R_3,\R_\tau\bigr] &= \R_z,      &\quad
\bigl[\R_3,\R_z\bigr]     &= 2\R_4.
\end{alignedat}
\end{equation}
In particular, \([\R_\tau,\R_z]=0\), and \(\R_4\) commutes with each of the six fields.
\end{proposition}
\begin{proof}
The brackets with $\R_2,\R_3$ follow by differentiating coefficients. In particular
\[
\begin{aligned}
[\partial_{t_1},\R_\tau]
 &=2a\partial_a+3b\partial_b+c\partial_c+2\partial_d+2t_1\partial_{t_1}+4t_2\partial_{t_2}
 =-\R_1.
\end{aligned}
\]
Let \(D=[\R_\tau,\R_z]\) and \(B=6a^2-t_2/2\). To check \(D=0\), it suffices
to evaluate it on the six coordinates. Formula~\eqref{eq:fields} gives
\[
\begin{aligned}
D(a)&=A_b-\R_zA_a=0,\qquad D(d)=2\R_\tau c-\R_z(c^2-a+2t_1)=0,\\
D(c)&=A_a+t_1^2-t_2/12-\R_z\bigl((a+t_1)c-b/2\bigr)=0.
\end{aligned}
\]
It also gives \(D(t_1)=D(t_2)=0\), since \(\R_z t_3=0\), and
\(D(b)=\R_\tau B-\R_z A_b=-12a^3+3at_2+3b^2+3t_3=0\)
by the equation of \(P\).
Finally, \(\R_1\) is minus the Euler field for the weights \((2,3,1,0,2,4)\)
of \((a,b,c,d,t_1,t_2)\), minus \(2\partial_d\).
The degrees of \(\R_\tau,\R_z,\R_2,\R_3\) are \(2,1,-2,-1\), respectively.
This gives all brackets with \(\R_1\); \(\R_4=\partial_d\) is central because the coefficients
in \eqref{eq:fields} are independent of \(d\). The remaining brackets vanish.
\end{proof}
\begin{corollary}\label{cor:flatness-lie-algebra}
The connection \eqref{eq:connection} is flat, and the complex Lie algebra spanned by the six vector fields is isomorphic to $\mathfrak{sl}_2\ltimes\mathfrak h_3$.
\end{corollary}
\begin{proof}
The matrices \eqref{eq:matrices} satisfy \eqref{eq:anti} for all the brackets
in \eqref{eq:brackets}, including the zero brackets. The curvature therefore
vanishes on a basis at every point. The triple \((-\R_1,\R_\tau,-\R_2)\)
satisfies the standard \(\mathfrak{sl}_2\) relations. The other three fields
span a Heisenberg ideal: \([\R_3,\R_z]=2\R_4\), and \(\R_4\) commutes with
all six fields. On the ordered basis \((\R_z,-\R_3)\) the triple acts by
the standard two dimensional representation and fixes \(\R_4\).
\end{proof}

\begin{remark}
The derivations also satisfy
\begin{equation}\label{eq:discriminant}
\begin{alignedat}{3}
\R_\tau t_3 &= 6t_1t_3-\frac{t_2^2}{3}, &\qquad
\R_\tau\Delta &= 12t_1\Delta, &\qquad
\R_z t_3 &= 0,\\[2pt]
\R_z\Delta &= 0, &\qquad
\R_1\Delta &= -12\Delta, &\qquad
\R_i\Delta &= 0\quad(i=2,3,4).
\end{alignedat}
\end{equation}
In particular, $\R_\tau(\Delta^{-1})=-12t_1\Delta^{-1}$. All six fields therefore act on $A_0[d]$ and, after $\partial_d\mapsto-h\partial_h$, on $A_0[h,h^{-1}]$. Together with \eqref{eq:det}, this proves the algebraic assertions of
Theorem~\ref{thm:main}.
\end{remark}

\section{Generalized periods}\label{sec:periods}


We compute the normalized periods, compare the period map with the
index zero case, and relate its connection to the algebraic matrix of
Section~\ref{sec:connection}.

Fix a symplectic marking of the elliptic curve. The relative path
\(\delta_0\) runs from \(\partial_O\) to \(\partial_P\) in the sense of
tangential basepoints. Choose lifts $\delta_1,\delta_2$ of the marked elliptic cycles to $X\setminus\{O,P\}$ compatible with this path and the
intersection convention below. Let $\delta_3$ be a positively
oriented small loop about $P$. The four paths and cycles
index the rows of the period matrix, while the ordered frame
\((\alpha_0,\alpha_1,\alpha_2,\alpha_3)\) indexes its columns. We choose logarithm branches on simply connected marked charts. Changes of the lift of \(P\), the paths, or these branches act through the discrete group \(\Gamma\) defined below. Throughout we use the tangents
\(\partial_O,\partial_P\) fixed in the calculation of the Gauss-Manin connection.

\subsection{The period domain}

We describe the normalized period matrices and the group acting on them by changes of marking.
The graph class has relative period \(-1\) and vanishes on closed cycles.
A differential with residue \(1\) at \(P\) has period \(2\pi i\) around the positively oriented loop \(\delta_3\).
Multiplication of its column by \(k^2=(-2\pi i)^{-1}\) makes this entry \(-1\).
The bilinear relations and the weight filtrations lead to the following
period domain $\Pi$:
\begin{equation}\label{eq:Pi}
 \mathcal P=
 \begin{pmatrix}
 -1 & z_1 & z_2 & w \\
  0 & x_1 & x_2 & x_2z_1-x_1z_2 \\
  0 & x_3 & x_4 & x_4z_1-x_3z_2 \\
  0 &   0 &   0 & -1
 \end{pmatrix},
\end{equation}
where \(x_1x_4-x_2x_3=1\), \(\operatorname{Im}(x_1\overline{x_3})>0\),
and \(z_1\notin\Z x_1+\Z x_3\). The strict inequality forces
\(x_3\ne0\), and then
\(\operatorname{Im}(x_1/x_3)=\operatorname{Im}(x_1\overline{x_3})/|x_3|^2>0\).
Since \(x_2=(x_1x_4-1)/x_3\), the six entries
\(x_1,x_3,x_4,z_1,z_2,w\) are coordinates on \(\Pi\).
Direct multiplication yields $\mathcal P^{\tr}\Phi\mathcal P=\Phi$.
This identity determines the
second and third entries of the last column in \eqref{eq:Pi}; it imposes
no condition on the upper right entry \(w\in\C\).
Let
\begin{equation}\label{eq:Gamma}
 \Gamma=\left\{
 \begin{pmatrix}
 1 & \lambda & \mu & \kappa \\
 0 & A & B & \mu' \\
 0 & C & D & -\lambda' \\
 0 & 0 & 0 & 1
 \end{pmatrix}
 \;\middle|\;
 \begin{gathered}
 AD-BC=1,\quad A,B,C,D,\lambda,\mu,\kappa\in\Z,\\
 \lambda'=\lambda D-\mu C,\quad
 \mu'=-\lambda B+\mu A
 \end{gathered}
 \right\}.
\end{equation}
\begin{proposition}\label{prop:Gamma-action}
Left multiplication defines a free, properly discontinuous holomorphic
action of $\Gamma$ on $\Pi$. Thus $\Gamma\backslash\Pi$ is a complex
manifold, and $\Pi\to\Gamma\backslash\Pi$ is a covering map.
\end{proposition}
\begin{proof}
Write \(\mathcal P'=\gamma\mathcal P\), and denote its entries as in
\eqref{eq:Pi}, with primes. Its first column and last row are those of
\(\mathcal P\). The middle \(2\times2\) matrix is multiplied on the left
by \(S=\left(\begin{smallmatrix}A&B\\C&D\end{smallmatrix}\right)\in\SL_2(\Z)\).
Thus its determinant is one, \(x'_1/x'_3=S\cdot(x_1/x_3)\in\HH\), and
\(\Z x'_1+\Z x'_3=\Z x_1+\Z x_3\). Since
\(z'_1=z_1+\lambda x_1+\mu x_3\), it is still outside this lattice.
The relations in \eqref{eq:Gamma} give \(\gamma^{\tr}\Phi\gamma=\Phi\).
Hence \((\mathcal P')^{\tr}\Phi\mathcal P'=\Phi\); pairing the last column
with the second and third columns gives its two middle entries as
\(x'_2z'_1-x'_1z'_2\) and \(x'_4z'_1-x'_3z'_2\), respectively.
This proves \(\gamma\mathcal P\in\Pi\). Invertibility of \(\mathcal P\)
gives freeness. For compact \(K\subset\Pi\), if \(\gamma K\cap K\ne\varnothing\),
then \(\gamma\in KK^{-1}\), a compact set containing only finitely many integral matrices.
\end{proof}

\begin{remark}
The matrix \(S=\left(\begin{smallmatrix}A&B\\ C&D\end{smallmatrix}\right)\)
changes the oriented basis of \(H_1(X,\Z)\). The first row of \(\gamma\)
replaces the relative path by
\(\delta_0+\lambda\delta_1+\mu\delta_2+\kappa\delta_3\).
Thus \(\lambda\) and \(\mu\) change the lift of \(P\) by an elliptic period.
Preservation of the pairing determines the multiples \(\mu'\) and
\(-\lambda'\) of the residue loop added to the two elliptic cycles in
\eqref{eq:Gamma}. Adding \(\kappa\delta_3\) alone translates \(w\) by \(-\kappa\).
\end{remark}

The quotient
\[
 \mathcal H_\Gamma=\Gamma\backslash(\Pi\times\C^4),
 \qquad
 \gamma\cdot(\mathcal P,v)=(\gamma\mathcal P,\gamma v),
\]
is a holomorphic vector bundle on \(\Gamma\backslash\Pi\), with the flat
connection induced by the trivial connection on \(\Pi\times\C^4\). For
the standard basis \(e_0,\ldots,e_3\) of \(\C^4\), the sections
\(\mathcal P\mapsto[\mathcal P,\mathcal P e_j]\) form a global
holomorphic frame. In the column frame convention, the connection matrix
with respect to this frame is
\[
 \A_\Pi=(\mathcal P^{-1}\dd\mathcal P)^{\tr}.
\]
\begin{proposition}
For each \(i\in\{\tau,z,1,2,3,4\}\), there is a unique holomorphic vector
field \(\mathcal R_i\) on \(\Pi\) satisfying \(\A_\Pi(\mathcal R_i)=C_i\).
It descends to \(\Gamma\backslash\Pi\). In the coordinates of \eqref{eq:Pi},
these fields are
\begin{equation}\label{eq:period-fields}
\begin{aligned}
 \mathcal R_\tau
   &=-x_2\partial_{x_1}-x_4\partial_{x_3}-z_2\partial_{z_1},\\
 \mathcal R_1
   &=x_1\partial_{x_1}-x_2\partial_{x_2}
     +x_3\partial_{x_3}-x_4\partial_{x_4}
     +z_1\partial_{z_1}-z_2\partial_{z_2},\\
 \mathcal R_2
   &=x_1\partial_{x_2}+x_3\partial_{x_4}+z_1\partial_{z_2},\\
 \mathcal R_z&=\partial_{z_1}+z_2\partial_w,\\
 \mathcal R_3&=-\partial_{z_2}+z_1\partial_w,\\
 \mathcal R_4&=-\partial_w.
\end{aligned}
\end{equation}
\end{proposition}
\begin{proof}
The defining equation is equivalent to
\(\dd\mathcal P(\mathcal R_i)=\mathcal P C_i^{\tr}\).
The entries \(x_1,x_3,x_4,z_1,z_2,w\) are local coordinates, so their
derivatives determine the field. Multiplication by the matrices in
\eqref{eq:matrices} gives \eqref{eq:period-fields}; substitution in
\eqref{eq:Pi} verifies the remaining entries. Left multiplication by
\(\gamma\in\Gamma\) preserves \(\A_\Pi\), so uniqueness gives invariance.
\end{proof}
By \eqref{eq:normalized-connection}, the differential of the normalized
period map sends \(\R_\tau,\R_1,\R_2,\R_z,\R_3,\R_4\) to
\(\mathcal R_\tau,\mathcal R_1,\mathcal R_2,k\mathcal R_z,
k\mathcal R_3,k^2\mathcal R_4\), respectively.

\subsection{Explicit normalized periods}

We compute the periods of the frame using the addition and quasi periodicity formulas for the
Weierstrass functions \cite{Chandrasekharan,DLMF}. We use full periods and regularize
the integral of the differential of the third kind along \(\delta_0\) using the endpoint
tangents fixed in Section~\ref{sec:cohom}.

Uniformize $X=\C/\Lambda$ by
\[
 \Lambda=\Z\omega_1+\Z\omega_2,\quad
 \operatorname{Im}(\omega_1/\omega_2)>0,\qquad
 (x,y,\alpha_1)=(\wp_\Lambda(u),\wp'_\Lambda(u),\dd u).
\]
Let $v\notin\Lambda$ lift $P$ and write
$\zeta_v=\zeta_\Lambda(v)$. Define the quasi periods for the full periods \(\omega_j\) by
$\zeta_\Lambda(u+\omega_j)=\zeta_\Lambda(u)+\eta_j$. Orient the cycles so that
\begin{equation}\label{eq:Legendre}
 \eta_1\omega_2-\eta_2\omega_1=-2\pi i.
\end{equation}
\begin{lemma}
\label{lem:elliptic-relative-periods}
Let \(\delta_0\) be the relative path with lift from \(0\) to \(v\).
For \(j=1,2\), let \(\delta_j\) be a closed cycle in
\(X\setminus\{O,P\}\) whose image in \(H_1(X,\Z)\) corresponds to
\(\omega_j\in\Lambda\). Then
\begin{equation}\label{eq:rawperiods}
\begin{alignedat}{2}
 \int_{\delta_0}\alpha_1&=v,
 &\qquad\int_{\delta_j}\alpha_1&=\omega_j,\\
 \int_{\delta_0}\alpha_2&=t_1v-\zeta_v-c,
 &\int_{\delta_j}\alpha_2&=t_1\omega_j-\eta_j
 \qquad(j=1,2).
\end{alignedat}
\end{equation}
The formulas extend holomorphically through $a=0$.
\end{lemma}
\begin{proof}
The first formulas follow from $\alpha_1=\dd u$.  For $\alpha_2$, the
primitive of $x\dd u$ is $-\zeta_\Lambda(u)$.  On $a\ne0$, the representative
\eqref{eq:frame} subtracts $y/(2x)$ and adds
$(c-b/(2a))\alpha_0$.  At $O$ the expression
$-\zeta_\Lambda(u)-y/(2x)$ has zero finite part, whereas at $P$ it equals
$-\zeta_v-b/(2a)$.  Since $\int_{\delta_0}\alpha_0=-1$, the terms involving
$b/(2a)$ cancel.  This proves \eqref{eq:rawperiods}.  At $a=0$ the apparent
pole is auxiliary; the local \v Cech representative of
Section~\ref{sec:moduli} proves the asserted continuation.
\end{proof}

Normalize \(\sigma_\Lambda\) by \(\sigma_\Lambda'(0)=1\), fix logarithms
on a simply connected marked chart, and put
\(\chi=\log\sigma_\Lambda(-v)-\log\sigma_\Lambda(v)\in(2\Z+1)\pi i\).
The zeta addition formula gives the primitive
\begin{equation}\label{eq:sigma-primitive}
\eta_P=\bigl(\zeta_\Lambda(u-v)-\zeta_\Lambda(u)+\zeta_v\bigr)\dd u
=\dd\bigl(\log\sigma_\Lambda(u-v)-\log\sigma_\Lambda(u)+u\zeta_v\bigr).
\end{equation}
\begin{lemma}\label{lem:third-kind-periods}
\begin{enumerate}
\item[(i)]
For \(j=1,2\), let \(\delta_j\) be any class in
\(H_1(X\setminus\{O,P\},\Z)\) whose image in \(H_1(X,\Z)\) is
the marked class with period \(\omega_j\) for \(\alpha_1\). Then
\[
\int_{\delta_j}\alpha_3\equiv(c+\zeta_v)\omega_j-\eta_jv
\pmod{2\pi i\Z}.
\]
Equality holds for the lift on which
\(\log(\sigma_\Lambda(u-v)/\sigma_\Lambda(u))\) has increment \(-\eta_jv\).
\item[(ii)] For the tangential path \(\delta_0\) and the chosen branches,
\begin{equation}\label{eq:relative-third}
\Reg\int_{\delta_0}\alpha_3=-2\log\sigma_\Lambda(v)+(c+\zeta_v)v-d-\chi .
\end{equation}
\end{enumerate} 
\end{lemma}
\begin{proof}
The transformation formula for \(\sigma_\Lambda\) under lattice translation gives the increment in \textup{(i)} modulo
\(2\pi i\Z\). Adding a residue loop changes the period by \(2\pi i\),
so the specified lift gives equality.
For \textup{(ii)}, denote the primitive in \eqref{eq:sigma-primitive} by
\(J(u)\). The parameters \(u\) at \(O\) and \(u-v\) at \(P\) are normalized by the endpoint tangents. Hence
\[
\begin{aligned}
J(\eps_O)&=\log\sigma_\Lambda(-v)-\log\eps_O+o(1),\\
J(v+\eps_P)&=\log\eps_P-\log\sigma_\Lambda(v)+v\zeta_v+o(1).
\end{aligned}
\]
Taking the constant term as in \eqref{eq:reg}, using
$\log\sigma_\Lambda(-v)=\log\sigma_\Lambda(v)+\chi$, and adding the
contributions $cv$ and $-d$ proves \eqref{eq:relative-third}.
\end{proof}

Set \(c_{\mathrm{per}}=2\log(2\pi k)-\chi\) and
\(w_{\mathrm{reg}}=\Reg\int_{\delta_0}\alpha_3-c_{\mathrm{per}}\),
with the branches of Lemma~\ref{lem:third-kind-periods}.
\begin{proposition}\label{prop:normalized-period-matrix}
After replacing \(\alpha_3\) by \(\alpha_3+c_{\mathrm{per}}\alpha_0\),
the period matrix and its normalization are
\begin{equation}\label{eq:Q}
 \begin{gathered}
 Q=\begin{pmatrix}
 -1 & v & t_1v-\zeta_v-c & w_{\mathrm{reg}} \\
  0 & \omega_1 & t_1\omega_1-\eta_1
    & -\eta_1v+(c+\zeta_v)\omega_1 \\
  0 & \omega_2 & t_1\omega_2-\eta_2
    & -\eta_2v+(c+\zeta_v)\omega_2 \\
  0 & 0 & 0 & 2\pi i
 \end{pmatrix},\\
 \mathcal P=QD_0,\qquad D_0=\operatorname{diag}(1,k,k,k^2).
 \end{gathered}
\end{equation}
They satisfy
\begin{equation}\label{eq:period-cup-identity}
 \mathcal P\in\Pi,\qquad Q^{\tr}\Phi Q=-2\pi i\Phi.
\end{equation}
\end{proposition}
\begin{proof}
Since $(2\pi k)^2=2\pi i$, one has $e^{-c_{\mathrm{per}}}=k^2$.  The constant frame
change subtracts $c_{\mathrm{per}}$ from the relative period and leaves
all closed periods and the connection unchanged. The equality
$e^{-c_{\mathrm{per}}}=k^2$ accounts for the fixed factor $k^2$
in the effective tensor~\eqref{eq:effective-tensor}.  Lemmas~\ref{lem:elliptic-relative-periods} and
\ref{lem:third-kind-periods} give the entries of \eqref{eq:Q}.  In the
coordinates of \eqref{eq:Pi}, its elliptic block is
\[
 \begin{pmatrix}x_1&x_2\\x_3&x_4\end{pmatrix}
 =k\begin{pmatrix}
    \omega_1&t_1\omega_1-\eta_1\\
    \omega_2&t_1\omega_2-\eta_2
   \end{pmatrix}.
\]
The other independent entries are $z_1=kv$,
$z_2=k(t_1v-\zeta_v-c)$, and $w=k^2w_{\mathrm{reg}}$.
Equation \eqref{eq:Legendre} implies that the elliptic block has determinant
$1$; substitution gives
$x_{2j}z_1-x_{2j-1}z_2
=k^2[-\eta_jv+(c+\zeta_v)\omega_j]$ for $j=1,2$. Finally
$k^2(2\pi i)=-1$.  Hence $\mathcal P\in\Pi$, and direct multiplication gives
$Q^{\tr}\Phi Q=-2\pi i\Phi$.
\end{proof}

\begin{corollary}\label{prop:cup-period-compatibility}
Suppose the endpoint tangents are \(v_O=\lambda\partial_O\) and \(v_P=\mu\partial_P\), and write
\(\alpha_3=\eta_P+c\alpha_1+d_T\alpha_0\). On a logarithmic chart,
the normalized period matrix is obtained from \eqref{eq:Q} by setting
\(d=d_T-\log\lambda-\log\mu\).
\end{corollary}
\begin{proof}
Proposition~\ref{prop:cup-tangent-calibration} gives the substitution.
Changing a logarithm by $2\pi i$ gives a central element of $\Gamma$.
\end{proof}

\subsection{The normalized period map}\label{sec:inverse}

Let \(\Pi_0\) and \(\Gamma_0\) be obtained from \(\Pi\) and \(\Gamma\)
by deleting the last row and column of each matrix. For a point of
\(\T_0^{\mathrm{an}}\), integration of \((\alpha_0,\alpha_1,\alpha_2)\)
along \((\delta_0,\delta_1,\delta_2)\), followed by multiplication on the
right by \(\operatorname{diag}(1,k,k)\), defines the period map for index
zero, \(\operatorname{Per}_0:\T_0^{\mathrm{an}}\to\Gamma_0\backslash\Pi_0\).
The fourth column of \eqref{eq:Q} contains the additional periods needed
on \(\T^{\times,\mathrm{an}}\). Choose a marking there, a lift of \(P\),
and a local branch of \(\log h\). Set \(d=-\log h\). Integrating
the resulting de Rham frame with the normalized endpoint tangents gives the matrix
\(\mathcal P=QD_0\) of \eqref{eq:Q}, hence a point of \(\Pi\). The
changes of marking and of logarithm described after \eqref{eq:Gamma}
act by left multiplication by elements of \(\Gamma\). Therefore the
construction descends to a natural holomorphic map
\[
 \operatorname{Per}:\T^{\times,\mathrm{an}}
 \longrightarrow \Gamma\backslash\Pi.
\]
A marked lift of a point \(x\) specifies an actual period matrix
\(\mathcal P\in\Pi\) representing \(\operatorname{Per}(x)\).
On the covering of such pairs \((x,\mathcal P)\), the marked period map
is \((x,\mathcal P)\mapsto\mathcal P\). It is \(\Gamma\) equivariant:
the action of \(\gamma\) sends \((x,\mathcal P)\) to
\((x,\gamma\mathcal P)\).
For the period map of an elliptic curve with its origin and a de Rham basis, obtained by integrating over a symplectic basis of \(H_1(X,\Z)\), see \cite[Sections~11.15--11.16]{MovBook}.

\begin{theorem}\label{thm:inverse}
The normalized period map induces a biholomorphism
\[
 \operatorname{Per}:\T^{\times,\mathrm{an}}
       \xrightarrow{\ \sim\ }\Gamma\backslash\Pi.
\]
\end{theorem}
\begin{proof}
By \eqref{eq:Gamma}, \(\Gamma_0\cong\SL_2(\Z)\ltimes\Z^2\), and
$1\to\Z\to\Gamma\to\Gamma_0\to1$ is the central extension in which
\(\kappa\in\Z\) translates
\(w\) by \(-\kappa\). In the commutative square below, \(p\) forgets the effective tensor and
\(q\) deletes the last row and column of a period matrix representative:
\[
\begin{CD}
\T^{\times,\mathrm{an}} @>{\operatorname{Per}}>>
                    \Gamma\backslash\Pi\\
@V{p}VV @VV{q}V\\
\T_0^{\mathrm{an}} @>{\operatorname{Per}_0}>>
                    \Gamma_0\backslash\Pi_0.
\end{CD}
\]
The bottom map is a biholomorphism by \cite[Section~8.3, Proposition~8.1]{CMVL}.
Indeed, if \(s=(2\pi i)^{-1/2}\) is the factor used there, their column
factors in genus one are \((1,s,-s)\), whereas ours are \((1,k,k)\).
The two matrices differ by right multiplication by
\(\operatorname{diag}(1,k/s,-k/s)\). Since \((k/s)^2=-1\), this
diagonal change preserves the determinant, the positivity and the
excluded lattice condition, and commutes with the left action of \(\Gamma_0\).

Fix \(x\in\T_0^{\mathrm{an}}\) and choose a matrix \(M\in\Pi_0\) representing \(\operatorname{Per}_0(x)\). If two matrices above \(M\) differ by \(\gamma\in\Gamma\), its image \(\gamma_0\in\Gamma_0\) satisfies \(\gamma_0M=M\). Since \(M\) is invertible, \(\gamma_0=I\). By \eqref{eq:Gamma}, the kernel of \(\Gamma\to\Gamma_0\) consists of the matrices \(I+nE_{03}\), which act by \(w\mapsto w-n\), \(n\in\Z\). Thus the two fibers have coordinates \(h\in\C^\times\) and \([w]\in\C/\Z\), respectively.

In a local marking, let \(w_*(x)\) be the value of \(w\) at \(h=1\).
Replacing \(1\) by \(h\) adds \(-\log(h)\alpha_0\) to \(\alpha_3\).
The relative period of \(\alpha_0\) is \(-1\), and the last column is
scaled by \(k^2\); hence \([w-w_*(x)]=[k^2\log h]\).
Since \(2\pi i k^2=-1\), the inverse on the fiber is
\(h=\exp((w-w_*(x))/k^2)\). This formula is holomorphic in the base
and fiber coordinates. By \eqref{eq:Gamma}, changing the marking adds
the same quantity to \(w\) and \(w_*(x)\) modulo \(\Z\).
The local inverses therefore agree and, together with
\(\operatorname{Per}_0^{-1}\), give a holomorphic inverse to \(\operatorname{Per}\).
\end{proof}

\subsection{The classical \texorpdfstring{$\tau$}{tau} map and geometric comparison}\label{sec:comparison}
We define the classical locus in \(\Pi\) and its lift to \(\T^{\times,\mathrm{an}}\). The period equation then relates differentiation along this locus to the vector fields of Section~\ref{sec:vector-fields}.
\begin{definition}\label{def:enhanced-t-map}
The \emph{$\tau$-period map} is
$\boldsymbol\tau:\mathcal H\to\Pi$, $(\tau,z)\mapsto\mathcal P(\tau,z)$,
where
\begin{equation}\label{eq:tau-locus}
 \mathcal P(\tau,z)=
 \begin{pmatrix}
 -1 & z & 0 & 0 \\
  0 & \tau & -1 & -z \\
  0 & 1 & 0 & 0 \\
  0 & 0 & 0 & -1
 \end{pmatrix}.
\end{equation}
The associated \emph{$\mathsf t$-map} is
$\mathsf t=\operatorname{Per}^{-1}\circ[\boldsymbol\tau]:
\mathcal H\to\T^{\times,\mathrm{an}}$,
where $\operatorname{Per}$ is the normalized period isomorphism of Theorem~\ref{thm:inverse}.
\end{definition}
The matrix \eqref{eq:tau-locus} belongs to $\Pi$. Following
\cite[Sections~8.4 and~9]{CMVL}, we restrict the period connection
to the image of $\boldsymbol\tau$.

\begin{lemma}\label{lem:slice-connection}
On \(\mathcal H\), let \(\mathcal P=\boldsymbol\tau(\tau,z)\) and
\(Q=\mathcal P D_0^{-1}\). Then
\[
(\mathcal P^{-1}\dd\mathcal P)^{\tr}=C_\tau\,\dd\tau+C_z\,\dd z,
\qquad
(Q^{-1}\dd Q)^{\tr}=C_\tau\,\dd\tau+k^{-1}C_z\,\dd z .
\]
\end{lemma}
\begin{proof}
Differentiating \eqref{eq:tau-locus} gives
\(\partial_\tau\mathcal P=\mathcal P C_\tau^{\tr}\) and
\(\partial_z\mathcal P=\mathcal P C_z^{\tr}\), by \eqref{eq:matrices}.
This gives the first identity. For the second, conjugate by \(D_0^{-1}\):
\(D_0^{-1}C_\tau D_0=C_\tau\) and \(D_0^{-1}C_zD_0=k^{-1}C_z\).
\end{proof}

\begin{proposition}\label{prop:analytic-connection}
On a simply connected open subset of \(\T^{\times,\mathrm{an}}\),
choose locally constant path and cycle classes \(\delta_0,\ldots,\delta_3\)
as above and a holomorphic function \(d\) with \(e^{-d}=h\). For the
corresponding period matrices \(Q\) and \(\mathcal P=QD_0\),
\begin{equation}\label{eq:normalized-connection}
 (Q^{-1}\dd Q)^{\tr}=\A,\qquad
 (\mathcal P^{-1}\dd\mathcal P)^{\tr}=D_0\A D_0^{-1}.
\end{equation}
\end{proposition}
\begin{proof}
Choose the path and cycle classes \(\delta_i\) as locally constant sections
of the integral homology local system on a simply connected marked chart.
The columns of \(Q\) are the periods of the frame in
Proposition~\ref{prop:normalized-period-matrix}. Differentiation of a
period over a locally constant class pairs that class with the covariant
derivative of the de Rham class. For \(\delta_0\), the endpoint terms
are those computed with the normalized tangents in Section~\ref{sec:connection}.
With cycles indexing rows and forms indexing columns, this gives
\[
 \dd Q=Q\A^{\tr}.
\]
Here $\A$ is the connection matrix calculated in
Section~\ref{sec:connection}. The constant frame change does not alter
that matrix: $\nabla\alpha_0=0$, and the last column of
\eqref{eq:connection} vanishes. Since $Q$ is invertible by
\eqref{eq:period-cup-identity}, this proves the first assertion.
The matrix $D_0$ is constant, so $\mathcal P=QD_0$ gives
\eqref{eq:normalized-connection} by conjugation.
\end{proof}

\begin{corollary}\label{prop:geometric}
The period isomorphism identifies the flat bundles
\[
(\mathcal H^\times,\nabla^\times)^{\mathrm{an}}
\xrightarrow{\ \sim\ }\operatorname{Per}^*\mathcal H_\Gamma,\qquad
\epsilon_j\longmapsto[\mathcal P,\mathcal P D_0^{-1}e_j].
\]
On a marked logarithmic chart this is the tangential period comparison
in the shifted frame of Proposition~\ref{prop:normalized-period-matrix}.
\end{corollary}
\begin{proof}
The displayed sections are independent of the lift $\mathcal P$:
both entries transform by $\gamma$. Their connection matrix is
$(Q^{-1}\dd Q)^{\tr}=\A^\times$ by
Proposition~\ref{prop:analytic-connection}. Moreover,
$Q^{\tr}(k^2\Phi)Q=\Phi$, so the constant pairing $k^2\Phi$ on the
Betti fiber induces the prescribed de Rham pairing. The entries of $Q$
are the compact, relative, and tangentially regularized periods of
the logarithmic frame after the constant frame change.

Replacing \(d\) by \(d+2\pi i n\) adds \(2\pi i n\alpha_0\) to
\(\alpha_3\). Since \(\alpha_0\) has relative period \(-1\) and
vanishing periods on closed cycles, only the upper right entry of \(Q\)
changes: it decreases by \(2\pi i n\). By \eqref{eq:Q}, the new matrix
is \((I-nE_{03})Q\). Since \(I-nE_{03}\in\Gamma\), the local
period sections obtained from the two branches agree in the quotient bundle.
For arbitrary endpoint tangents,
Proposition~\ref{prop:cup-tangent-calibration} replaces $d$ by
$d_T-\log\lambda-\log\mu$, whose differential is $-\dd h/h$.
\end{proof}

\section{Quasi Jacobi forms of integral index}\label{sec:Jacobi}
We determine the pullbacks of the regular coordinates on \(\T^\times\)
along \(\mathsf t\). As in the case of index zero \cite{CMVL}, changes
of marking give the transformation laws, and the Gauss-Manin
connection gives the differential equations. Theorem~\ref{prop:coordinate-algebra} characterizes the homogeneous
coordinate polynomials whose pullbacks are quasi Jacobi forms in
Oberdieck's sense. Proposition~\ref{prop:weak-ring} gives generators
and a defining relation for the ring of weak Jacobi forms of integral index and trivial multiplier.
For the Weierstrass \(\wp\) function, Jacobi's theta functions, and the Jacobi transformation laws, see \cite[Sections~2.4, 2.15, and~14.3]{MovBook}, respectively.

\subsection{Analytic preliminaries}\label{sec:analytic-Jacobi}
We fix the conventions for theta functions, the eta function, Weierstrass functions, and Jacobi forms. The heat equation supplies the identities needed to identify the Gauss-Manin vector fields.

Put $q=e^{2\pi i\tau}$ and $\zeta=e^{2\pi iz}$. We normalize the odd theta function and the eta function by
\[
 \begin{aligned}
 \theta_1(\tau,z)
   =-i\sum_{n\in\Z}(-1)^n
     e^{\pi i(n+1/2)^2\tau+2\pi i(n+1/2)z}, \qquad
 \eta(\tau)=q^{1/24}\prod_{n\ge1}(1-q^n).
 \end{aligned}
\]
Thus $\theta_1'(\tau,0)=2\pi\eta(\tau)^3$. Write $F=\theta_1/\eta^3$ and $J_1=\partial_z\log F$. The Jacobi product gives
\begin{equation}\label{eq:F-product}
 F(\tau,z)=2\sin(\pi z)
 \prod_{n\ge1}\frac{(1-\zeta q^n)(1-\zeta^{-1}q^n)}{(1-q^n)^2}.
\end{equation}
Hence $F$ is holomorphic on $\HH\times\C$, with simple zeros precisely at $z\in\Z\tau+\Z$. The function $J_1$ is meromorphic with a simple pole of residue $1$ at each such point. For $\lambda,\mu\in\Z$, the product gives
\begin{equation}\label{eq:F-shift}
 F(\tau,z+\lambda\tau+\mu)
 =(-1)^{\lambda+\mu}
 e^{-\pi i(\lambda^2\tau+2\lambda z)}F(\tau,z).
\end{equation}
We write $\wp(\tau,z)$ for the Weierstrass function of $\Z\tau+\Z$ and $\wp'(\tau,z)$ for its $z$ derivative.
The normalized Eisenstein series are
$E_{2r}=1-\frac{4r}{B_{2r}}\sum_{n\ge1}n^{2r-1}q^n/(1-q^n)$ for
$r=1,2,3$, where $B_{2r}$ denotes the corresponding Bernoulli number.

The theta series and eta product converge locally uniformly. The theta heat equation and eta product give
\[
 \partial_\tau\theta_1=\frac{1}{4\pi i}\partial_z^2\theta_1,
 \qquad \partial_\tau\log\eta=\frac{\pi i}{12}E_2.
\]
On any simply connected subset avoiding the lattice divisor, set $L=\log F$. Then
\begin{equation}\label{eq:heat}
 L_\tau=\frac{L_{zz}+L_z^2}{4\pi i}-\frac{\pi i}{4}E_2,
 \qquad L_{zz}=-\wp(\tau,z)-\frac{\pi^2}{3}E_2.
\end{equation}
Indeed, $L_{zz}+\wp$ is elliptic without poles, hence independent of $z$. Comparing the Laurent expansions at the origin in \eqref{eq:F-product} gives its constant as $-\pi^2E_2/3$. In particular, $J_1$ satisfies the analytic differential identities
\begin{equation}\label{eq:J1-analytic}
 \partial_zJ_1=-\wp-\frac{\pi^2}{3}E_2,
 \qquad 4\pi i\,\partial_\tau J_1
 =\partial_z\!\left(\partial_zJ_1+J_1^2\right).
\end{equation}
The second equation follows by differentiating \eqref{eq:heat} in $z$; it is independent of the choice of the local logarithm.

We use the following conventions for classical Jacobi forms. If
$\gamma=\begin{psmallmatrix}A&B\\C&D\end{psmallmatrix}\in\SL_2(\Z)$
and $(\lambda,\mu)\in\Z^2$, set
\begin{equation}\label{eq:Jacobi-action}
 \begin{aligned}
 j=C\tau+D,  \qquad v=z+\lambda\tau+\mu, \qquad 
 \tau'=\frac{A\tau+B}{j},  \qquad z'=\frac vj, \qquad
 Q_J=\frac{Cv^2}{j}-\lambda^2\tau-2\lambda z.
 \end{aligned}
\end{equation}
\begin{definition}\label{def:weak-Jacobi}
A meromorphic function \(f\) on \(\HH\times\C\) is a \emph{meromorphic
Jacobi function} of weight \(\ell\in\Z\) and index \(m\in\Z\), with
trivial multiplier, if every transformation \eqref{eq:Jacobi-action} satisfies
\[
f(\tau',z')=j^\ell e^{2\pi imQ_J}f(\tau,z).
\]
For \(m\ge0\), let \(\widetilde J_{\ell,m}\) denote the space of
\emph{weak Jacobi forms}: its elements are holomorphic functions on
\(\HH\times\C\) satisfying this law and admitting a locally uniformly
convergent Fourier expansion
\(f=\sum_{n\ge0,\,r\in\Z}c(n,r)q^n\zeta^r\) at the cusp.
The space \(J_{\ell,m}\) of \emph{holomorphic Jacobi forms} is the subspace
of \(\widetilde J_{\ell,m}\) defined by \(c(n,r)=0\) for \(r^2>4mn\).
\end{definition}

Following \cite[Appendix~B]{Oberdieck}, introduce the bigraded analytic algebra
\begin{equation}\label{eq:Oberdieck-algebra}
 \mathcal V=\C[F,E_2,E_4,J_1,\wp,\wp']
   =\bigoplus_{m\in\frac12\Z_{\ge0}}
      \bigoplus_{\ell\in\Z}\mathcal V_{\ell,m},
\end{equation}
where the weight index assignments are
\[
\begin{array}{c|rrrrrr}
 &F&E_2&E_4&J_1&\wp&\wp'\\ \hline
\mathrm{weight}&-1&2&4&1&2&3\\
\mathrm{index}&\frac12&0&0&0&0&0.
\end{array}
\]
Multiplying any generator by a fixed nonzero constant does not change this bigraded algebra; hence this definition is compatible with the $k$ normalization used in \eqref{eq:theta-map}.

\begin{definition}[Oberdieck]\label{def:quasi-Jacobi}
The space \(\operatorname{QJac}_{\ell,m}\) of quasi Jacobi forms of
weight \(\ell\) and index \(m\) consists of those \(f\in\mathcal V_{\ell,m}\)
that are holomorphic at \(z=0\) for generic \(\tau\in\HH\).
Set \(\operatorname{QJac}=\bigoplus_{\ell,m}\operatorname{QJac}_{\ell,m}\).
\end{definition}

\begin{remark}\label{rem:weak-quasi}
In Definition~\ref{def:quasi-Jacobi}, holomorphy is tested at the
origin. Fix \(\tau\) and suppress it from the notation. Then \(J_1(u)=u^{-1}+O(u)\) and \(\wp(u)=u^{-2}+O(u^2)\),
so \(J_1^2-\wp\) is regular there. At \(z=\tau+u\), however,
\(J_1(\tau+u)=J_1(u)-2\pi i\) by \eqref{eq:F-shift}, whereas
\(\wp(\tau+u)=\wp(u)\); the principal part is therefore \(-4\pi i/u\).
For a function satisfying the exact Jacobi law, holomorphy at the origin
does imply holomorphy at all its lattice translates.
\end{remark}

\subsection{The \texorpdfstring{$\tau$}{tau}-map, the enhanced period map and automorphy}\label{sec:t-map}
In this subsection, we derive the transformation laws of the enhanced coordinates from
changes of marking before identifying their pullbacks on the
$(\tau,z)$ domain. The transformation laws, Gauss-Manin equations,
and local normalization together determine those pullbacks; automorphy
alone leaves holomorphic factors undetermined.

The maps used here are the $\tau$-period map $\boldsymbol\tau$ and the enhanced map $\mathsf t$ of Definition~\ref{def:enhanced-t-map}.  Write \(g\in\Gamma\) as in \eqref{eq:Gamma}, with middle
\(2\times2\) matrix \(\begin{psmallmatrix}A&B\\ C&D\end{psmallmatrix}\)
and integers \(\lambda,\mu,\kappa\). Set
\[
 R_g=\lambda-\frac{Cv}{j},\qquad S_g=-Q_J+\lambda\mu+\kappa,
\]
where \(j=C\tau+D\), \(v=z+\lambda\tau+\mu\), and \(Q_J\) is as in
\eqref{eq:Jacobi-action}. The change of marking sends \((\tau,z)\) to
\[
 (\tau',z')=\left(\frac{A\tau+B}{C\tau+D},
                   \frac{z+\lambda\tau+\mu}{C\tau+D}\right).
\]
The matrix \(G_g\) below expresses the columns of \(g\mathcal P(\tau,z)\)
in the frame given by the columns of \(\mathcal P(\tau',z')\):
\begin{equation}\label{eq:Jacobi-factor}
 \begin{gathered}
 g\,\mathcal P(\tau,z)=\mathcal P(\tau',z')G_g,\\
 G_g=
 \begin{pmatrix}
 1 & 0 & R_g & S_g \\
 0 & j & -C & jR_g \\
 0 & 0 & j^{-1} & 0 \\
 0 & 0 & 0 & 1
 \end{pmatrix}.
 \end{gathered}
\end{equation}
For instance, comparison of the first rows gives $R_g=\lambda-Cv/j$ and $S_g=R_gv+\lambda z+\kappa$; substitution of \eqref{eq:Jacobi-action} yields the displayed expression for $S_g$. The remaining entries follow using $AD-BC=1$.

Left multiplication by \(g\) changes only the marking.
Equation~\eqref{eq:Jacobi-factor} therefore expresses the frame at \((\tau,z)\)
in terms of the frame at \((\tau',z')\). Before the column normalization by
\(D_0\), the matrix of this change is \(D_0G_gD_0^{-1}\). In the notation
of \eqref{eq:G}, the frame parameters are
\((s,s',v,w)=(j,-C,R_g/k,S_g/k^2)\).
Applying the inverse change of frame to the coordinates gives
\begin{equation}\label{eq:automorphy}
\begin{alignedat}{3}
 a'&=j^2a, &\qquad b'&=j^3b,
 &\qquad c'&=jc+\frac{Cv-j\lambda}{k},\\
 t'_1&=j^2t_1+Cj, & t'_2&=j^4t_2,
 & t'_3&=j^6t_3,\\
 h'&=j^{-2}e^{2\pi iQ_J}h,
 & d'&=d+2\log j-2\pi iQ_J
           \pmod{2\pi i\Z}.
\end{alignedat}
\end{equation}
Here a prime denotes evaluation at $(\tau',z')$, and the $d$ formula is on a chosen logarithmic chart. The identity for $t_3$ follows from $t_3=4a^3-t_2a-b^2$; in particular $\Delta'=j^{12}\Delta$. The central integer $\lambda\mu+\kappa$ drops out of the transformation of $h$, because $k^{-2}=-2\pi i$.

\begin{proposition}\label{prop:t-map-coordinates}
The pullbacks of the enhanced coordinates along $\mathsf t$ are given by
\eqref{eq:theta-map}.
On a logarithmic chart with compatible branches,
$\mathsf t^*d=-2\log F$.
\end{proposition}
\begin{proof}
We follow the method of \cite[Section~9.3]{CMVL}. The periods of
\(\alpha_1\) in \eqref{eq:Q} and \eqref{eq:tau-locus} give
\(\Lambda=k^{-1}(\Z\tau+\Z)\) and \(v=z/k\).
For fixed \(\tau\), the function \(\widetilde a=\mathsf t^*a\) is
elliptic in \(z\) by \eqref{eq:automorphy} and holomorphic off the lattice.
In the local coordinate \(u\) at \(O\) with \(\dd u=\alpha_1\) and
\(u(O)=0\), one has \(x(u)=u^{-2}+O(u^2)\). Since \(u(P)=z/k\) near
\(z=0\), the difference \(\widetilde a-k^2\wp(\tau,z)\) is an elliptic
function without poles and vanishes at the origin. Thus
\(\widetilde a=k^2\wp(\tau,z)\).

By Proposition~\ref{prop:analytic-connection} and Lemma~\ref{lem:slice-connection},
\(\dd\mathsf t(k\partial_z)\) solves the defining equation
\(\A(V)=C_z\) for \(\R_z\). Thus \(\mathsf t^*b=k\partial_z\widetilde a=k^3\wp'\).
Since \(\R_z t_2=\R_z t_3=0\), the pullbacks of \(t_2,t_3\) depend only
on \(\tau\). Comparing the coefficients of \(\wp\) and the constant
terms in the two Weierstrass equations gives
$\mathsf t^*t_2=k^4g_2(\Z\tau+\Z)=k^{-4}E_4/12$ and
$\mathsf t^*t_3=k^6g_3(\Z\tau+\Z)$.

Let \(\eta_0=\pi^2E_2/3\) be the quasi period of the zeta function for
\(\Z\tau+\Z\) under translation by \(1\). The entry \(x_4=0\) in \eqref{eq:tau-locus} and \eqref{eq:Q} give $t_1\omega_2=\eta_2$.
Scaling the quasi period to $\Lambda$ yields
$\mathsf t^*t_1=k^2\eta_0=-k^{-2}E_2/12$.
Moreover, $k\partial_z(\mathsf t^*c)=\widetilde a+\mathsf t^*t_1$.
Since $k^2\pi^2E_2/3=\mathsf t^*t_1$, equation~\eqref{eq:J1-analytic}
shows that $\mathsf t^*c+kJ_1$ is independent of $z$.
For $\gamma=-I$, \eqref{eq:automorphy} makes $\mathsf t^*c$
odd under $z\mapsto-z$. The function $J_1$ is odd as well,
so the constant vanishes and $\mathsf t^*c=-kJ_1$.

The theta transformation formula \cite[\S20.7(viii)]{DLMF}, applied to
\(F=2\pi\theta_1/\theta_1'(\tau,0)\), and \eqref{eq:F-shift} show that
\(F^2\) has weight \(-2\), index \(1\), and trivial multiplier. Thus
\eqref{eq:automorphy} makes \(\mathsf t^*h/F^2\) invariant, but does
not determine its value. The equation \(\R_z h=-2ch\) gives
\(\partial_z\log(\mathsf t^*h)=2J_1\), so the quotient is a function
\(C(\tau)\) of \(\tau\) alone. To determine it, use the normalized
tangents of Definition~\ref{def:normalized-tangents}.
In Lemma~\ref{lem:third-kind-periods}, substitute
\(c_{\mathrm{per}}=2\log(2\pi k)-\chi\) and
\(w=z_2=0\) on the image of \(\boldsymbol\tau\).
This gives \(c+\zeta_\Lambda(v)=t_1v\) and
\[
 0=-2\log\sigma_\Lambda(v)+t_1v^2-d-2\log(2\pi k).
\]
Thus $h=(2\pi k)^2\sigma_\Lambda(v)^2e^{-t_1v^2}$.
Since $v=z/k$ and $\sigma_\Lambda(v)=v+O(v^5)$, we obtain
$\mathsf t^*h=(2\pi z)^2(1+O(z^2))$ as $z\to0$.
Since $F=2\pi z+O(z^3)$, it follows that $C(\tau)=1$.
On a compatible logarithmic chart, $d=-\log h=-2\log F$.
\end{proof}

Let
\begin{equation}\label{eq:algebra}
B=A_0[h,h^{-1}],\qquad
\wt(a,b,c,t_1,t_2,h)=(2,3,1,2,4,-2),\qquad \ind(h)=1,
\end{equation}
with index zero on the other generators. Set $B_{\mathrm{pol}}=
 \frac{\C[a,b,c,t_1,t_2,t_3,h]}
 {\bigl(b^2-4a^3+t_2a+t_3\bigr)}.$

\begin{theorem}\label{prop:coordinate-algebra}
\begin{itemize}
\item The image of \(B_{\mathrm{pol}}\) under \(\mathsf t^*\) is
\(\bigoplus_{\ell\in\Z,\,m\in\Z_{\ge0}}\mathcal V_{\ell,m}\).
For \(f\in B_{\mathrm{pol}}\) of bidegree \((\ell,m)\), its pullback
belongs to \(\operatorname{QJac}_{\ell,m}\) precisely when it is
holomorphic at \(z=0\) for generic \(\tau\).
\item For \(m\in\Z\) and \(p\in A_0\) homogeneous of weight
\(\ell+2m\), the transformation \eqref{eq:Jacobi-action} gives
\begin{equation}\label{eq:quasi-law}
(h^mp)'=j^\ell e^{2\pi imQ_J}h^m
\sum_{r,s\ge0}\frac{u^rv_0^s}{r!s!}\partial_c^r\partial_{t_1}^sp,
\qquad u=\frac{Cv-j\lambda}{kj},\quad v_0=\frac Cj .
\end{equation}
\end{itemize}
\end{theorem}
\begin{proof}
The identities for the generators follow from Proposition~\ref{prop:t-map-coordinates}.
In integral index, every monomial in $\mathcal V$ contains an even
power of $F$; also $E_6$ is a polynomial in $\wp,\wp',E_4$ by the
Weierstrass equation.  Thus the image is precisely the integral index part of $\mathcal V$.
The holomorphy criterion is Oberdieck's definition.  More explicitly,
the negative Laurent coefficients at $z=0$ are holomorphic functions
of $\tau$; if they vanish generically, the identity theorem makes
them vanish identically.  For \eqref{eq:quasi-law}, factor out the
homogeneous powers of $j$ in \eqref{eq:automorphy}, then Taylor expand
in the shifts of $c$ and $t_1$.  Since $\Delta$ is independent of
$c,t_1$, the expansion is finite even in $A_0$.
\end{proof}

Every bihomogeneous element of \(B\) independent of \(c,t_1\)
satisfies the exact Jacobi transformation law of its weight and index
after pullback along \(\mathsf t\). Passing to $B$ permits inverses of $F^2$ and the
modular discriminant; holomorphy along the lattice divisor and at the
cusp must therefore be checked separately.

\begin{proposition}\label{prop:index-quotient}
For $\ell,m\in\Z$, multiplication by $(\mathsf t^*h)^m$ is a bijection
between meromorphic functions on $\mathcal H$ satisfying the Jacobi
transformation law of weight $\ell+2m$ and index zero and those satisfying
the law of weight $\ell$ and index $m$, in both cases with trivial multiplier.
\end{proposition}
\begin{proof}
The pullback $\mathsf t^*h=F^2$ is nowhere zero on $\mathcal H$.  By \eqref{eq:automorphy}, multiplication by \((\mathsf t^*h)^m\) gives the asserted transformation law; division by the same factor gives the inverse.
\end{proof}
For a meromorphic Jacobi function \(f\) on \(\mathcal H\), the quotient
\(f/F^{2m}\) need not extend across the lattice divisor or be the
pullback of an element of \(A_0\).

\subsection{The ring of weak Jacobi forms}\label{sec:weak-ring}
We express the classical structure theorem for weak Jacobi forms
\cite[Theorems~9.3 and~9.4]{EZ} in enhanced coordinates. Near \(z=0\)
one has
\[
 \begin{alignedat}{2}
 a=k^2z^{-2}+O(z^2), \qquad b=-2k^3z^{-3}+O(z), \qquad
 c=-kz^{-1}+O(z),  \qquad h=4\pi^2z^2+O(z^4).
 \end{alignedat}
\]
Thus $bh$ has a simple pole, whereas $h$, $ah$ and $bh^2$
extend holomorphically across the lattice.

Put \(U=h,\ V=ah,\ W=bh^2\); the bidegrees of \((t_2,t_3,U,V,W)\) are
\((4,0)\), \((6,0)\), \((-2,1)\), \((0,1)\) and \((-1,2)\), respectively.
\begin{proposition}\label{prop:weak-ring}
For trivial multiplier, pullback along \(\mathsf t\) induces an
isomorphism of bigraded algebras
\begin{equation}\label{eq:weak-ring}
\frac{\C[t_2,t_3,U,V,W]}{(W^2-4V^3U+t_2VU^3+t_3U^4)}
\;\xrightarrow{\sim}\;
\bigoplus_{\ell\in\Z,\,m\in\Z_{\ge0}}\widetilde J_{\ell,m}.
\end{equation}
\end{proposition}
\begin{proof}
By Proposition~\ref{prop:t-map-coordinates}, the pullbacks of
\(t_2,t_3,U,V,W\) are nonzero constant multiples of
\(E_4,E_6,\phi_{-2,1},\phi_{0,1},\phi_{-1,2}\), respectively.
The even weight ring is polynomial in the first four forms, and the
odd weight forms constitute a free module of rank one over this ring,
generated by \(\phi_{-1,2}\)
\cite[Theorems~9.3 and~9.4]{EZ}. Thus, writing
\(S=\C[t_2,t_3,U,V]\) and identifying these generators with their
pullbacks, the full ring is \(S\oplus WS\). The Weierstrass equation gives
\[
W^2=b^2h^4=4V^3U-t_2VU^3-t_3U^4.
\]
The direct sum decomposition shows that there is no further relation.
For fixed index $m\ge1$, a basis over $\C[t_2,t_3]$ consists of
$U^{m-j}V^j$ for $0\le j\le m$ and
$WU^{m-2-j}V^j$ for $0\le j\le m-2$;
the second range is empty when $m=1$.
Dividing by $h^m$ gives the normal form in
Corollary~\ref{thm:intro-FJ} and its pole bound $2m$.
\end{proof}

\begin{remark}\label{rem:EZ-ring}
In the standard normalization,
\(\phi_{0,1}=3V/(\pi^2k^2)=\zeta+10+\zeta^{-1}+O(q)\).
The Fourier coefficient condition $r^2\le4mn$ selects the
holomorphic Jacobi subring.
\end{remark}
If \(f\) satisfies the Jacobi transformation law of index \(m<0\)
and is holomorphic on \(\HH\times\C\), then \(f=0\). Indeed, for
fixed \(\tau\), the function \(fF^{-2m}\) is holomorphic and
elliptic. It is constant and vanishes at \(z=0\).

\subsection{Gauss-Manin vector fields and analytic differential equations}\label{sec:theta-proof}
The period identity of Proposition~\ref{prop:analytic-connection}
already identifies the tangent directions along the image of
$\boldsymbol\tau$.
The calculation below verifies the resulting differential identities
directly from the classical functions and determines the action at
every integral index.

\begin{proposition}\label{prop:analytic-vector-fields}
Along the enhanced map $\mathsf t$ one has
\[
 \dd\mathsf t(k\partial_z)=\R_z\circ\mathsf t, \qquad
 \dd\mathsf t(\partial_\tau)=\R_\tau\circ\mathsf t.
\]
\end{proposition}
\begin{proof}
Put $D_z=k\partial_z$ and locally $d=-2L=-2\log F$. The Weierstrass identities
$(\wp')^2=4\wp^3-g_2\wp-g_3$, $\wp''=6\wp^2-g_2/2$,
$g_2=4\pi^4E_4/3$ and $g_3=8\pi^6E_6/27$ imply
\begin{equation}\label{eq:theta-z}
\begin{alignedat}{4}
D_za &= b, &\qquad
D_zb &= 6a^2-\tfrac{t_2}{2}, &\qquad
D_zc &= a+t_1, &\qquad
D_zd &= 2c,\\[2pt]
D_zt_1 &= 0, &\qquad
D_zt_2 &= 0, &\qquad
D_zh &= -2ch.
\end{alignedat}
\end{equation}
Here and below coordinate symbols mean their pullbacks under $\mathsf t$; the equation for $d$ is local and the equation for $h$ is global. For example, $D_zc=-k^2L_{zz}=a+t_1$ by \eqref{eq:heat}, and $D_zh=2kL_z h=-2ch$.

The heat equation gives
\[
 \partial_\tau d=-2L_\tau
 =k^2(L_{zz}+L_z^2)+\frac{\pi i}{2}E_2
 =c^2-a+2t_1.
\]
Differentiating in $z$ and using $D_zd=2c$ yields
$\partial_\tau c=(a+t_1)c-b/2$. 
Ramanujan’s differential identities
\[
\partial_\tau E_2=\frac{2\pi i}{12}(E_2^2-E_4),\qquad
\partial_\tau E_4=\frac{2\pi i}{3}(E_2E_4-E_6)
\]
give \(\partial_\tau t_1=t_1^2-t_2/12\) and
\(\partial_\tau t_2=4t_1t_2-6t_3\).

Finally $a=D_zc-t_1$, $b=D_za$ and $[D_z,\partial_\tau]=0$ give
\begin{equation}\label{eq:theta-tau}
\begin{aligned}
 \partial_\tau a=-2a^2+2at_1+bc+t_2/3,
 \partial_\tau b=(6a^2-t_2/2)c-3ab+3bt_1,
 \partial_\tau h=-(c^2-a+2t_1)h.
\end{aligned}
\end{equation}
These identities give both asserted equalities of vector fields
along \(\mathsf t\). The discriminant is nonzero on the image because
\[
 \Delta=-\frac{k^{-12}}{1728}(E_4^3-E_6^2)
        =-k^{-12}q\prod_{n\ge1}(1-q^n)^{24}.
\]
\end{proof}

These vector fields act algebraically for every integer index, including
on expressions that are not classical weak Jacobi forms.  Put \(A_0=\C[a,b,c,t_1,t_2,\Delta^{-1}]\), where
\(t_3=4a^3-t_2a-b^2\) and \(\Delta=27t_3^2-t_2^3\).
On \(B=A_0[h,h^{-1}]\), their action on \(h\) is
\begin{equation}\label{eq:h-fields}
 \begin{alignedat}{3}
 \R_z h&=-2ch,
 &\qquad \R_\tau h&=-(c^2-a+2t_1)h,\\
 \R_1h&=2h,
 & \R_2h&=\R_3h=0,
 &\qquad \R_4h&=-h.
 \end{alignedat}
\end{equation}
On \(h^m A_0\), the index operator \(\mathcal I=h\partial_h=-\R_4\)
acts by \(m\); on elements of weight \(\ell\), the weight operator
\(-\R_1\) acts by \(\ell\). Write \(\R_z^0,\R_\tau^0\) for the
restrictions of \(\R_z,\R_\tau\) to \(A_0\). In the basis \(h^m\) of
the free \(A_0\) module \(h^m A_0\), these vector fields act by
\begin{equation}\label{eq:index-operators}
 D_z^{(m)}=\R_z^0-2mc,\qquad
 D_\tau^{(m)}=\R_\tau^0-m(c^2-a+2t_1),
\end{equation}
where the final terms act by multiplication. Indeed,
\[
 \R_z(h^mp)=h^mD_z^{(m)}p,\qquad
 \R_\tau(h^mp)=h^mD_\tau^{(m)}p.
\]
These operators commute, since \([\R_z^0,\R_\tau^0]=0\) and
\[
 \R_z^0(c^2-a+2t_1)=2(a+t_1)c-b=2\R_\tau^0c.
\]
Moreover, \([\partial_c,D_z^{(m)}]=-2m\). On \(h^m A_0\) the central
additive frame parameter \(w\) acts by \(e^{-mw}\), so \(\R_4\) acts by
\(-m\). For the multiplicative parameter \(\lambda=e^{-w}\), the
character is \(\lambda\mapsto\lambda^m\); \(2\pi i\Z\) acts trivially
in the additive parameter. After restoring \(h^m\), the operators
correspond, under \(\mathsf t\), to \(k\partial_z\) and \(\partial_\tau\).
Both derivatives preserve the integral index part of \(\operatorname{QJac}\):
they preserve its ambient polynomial algebra, its index, and holomorphy at
\(z=0\). They need not preserve the exact Jacobi law, as
\(k\partial_z F^2=2kJ_1F^2\) shows.
\begin{remark}
For a modular form of weight \(\ell\), the Serre derivative is
\((2\pi i)^{-1}(\partial_\tau-\ell\,\mathsf t^*t_1)\), since
\(\mathsf t^*t_1=(2\pi i)E_2/12\).
For a Jacobi form of positive index \(m\), the classical corrected heat
operator is
\[
\frac{1}{2\pi i}\partial_\tau
-\frac{1}{4m}\left(\frac{1}{2\pi i}\partial_z\right)^2
-\frac{2\ell-1}{24}E_2.
\]
Indeed, in the theta decomposition it acts on each coefficient by the
Serre derivative of weight \(\ell-\tfrac12\).
Thus \(D_\tau^{(m)}\) represents ordinary differentiation in the basis \(h^m\); the heat operator sends Jacobi forms of weight \(\ell\) and index \(m\) to Jacobi forms of weight \(\ell+2\) and index \(m\). For the Jacobi-Serre derivative in even weight, see
\cite[Theorem~1]{OberdieckSerre}.
\end{remark}

\section{Degenerations at a single node and boundary values}\label{sec:degeneration-new}
We first recall how the period matrix of a genus two curve behaves when the curve acquires a single node.
Each singular fiber of the family \eqref{eq:explicit-quintic} has one
nonseparating node, and its normalization is an elliptic curve. The
smoothing parameter determines a tensor of the tangent lines at the two
preimages of the node. Using this tensor and the limiting de Rham frame,
we identify the regularized periods and limiting connection with those of
the enhanced elliptic curve and compute the boundary values of the genus
two generators in the chosen frame.
Comparing the smoothing parameter with the Siegel period parameter then expresses the
leading term of a Siegel modular form in the chosen smoothing parameter. We treat a separating node at the end.
For background on period maps, see \cite[Chapter~1]{CMSP}.

\subsection{Plumbing coordinates and normalized periods}


Let \(\pi:\mathcal C\to S\) be a semistable family of genus two curves over a smooth local
analytic base. Suppose that the singular fibers lie over a smooth divisor \(D=\{s=0\}\),
that each has just one node, and that \(\mathcal C\) is smooth at these nodes.
Fix \(x_0\in D\), write \(C_0=\mathcal C_{x_0}\), and denote the preimages of its node in the
normalization \(\widetilde C_0\) by \(O,P\). Choose local parameters \(u\) at \(O\) and
\(v\) at \(P\) in which the family has equation
\begin{equation}\label{eq:plumbing-section8}
uv=s.
\end{equation}
Over a sector in the punctured \(s\) disk, fix a branch of \(\log s\) and a homology marking.
The other base coordinates are held fixed when taking limits in \(s\).
\begin{lemma}\label{lem:plumbing-tensor}
A normal vector \(\xi\in N_{D/S,x_0}\) determines the tensor
\[
\upsilon^{\mathrm{sm}}_\xi
=\dd s(\xi)\,\partial_u|_O\otimes\partial_v|_P
\in T_O\widetilde C_0\otimes T_P\widetilde C_0,
\]
independently of the chosen plumbing coordinates. In a second chart \(u'v'=s'\), take \(u'\) to be a local parameter at
\(O\) and \(v'\) a local parameter at \(P\), as for \(u\) and \(v\),
respectively. Along \(D\), let \(A=(\dd u'/\dd u)|_O\) and
\(B=(\dd v'/\dd v)|_P\) be the derivatives on the two branches. With \(I_D\) the ideal sheaf of \(D\), we have
\begin{equation}\label{eq:plumbing-coordinate-change}
[s']=AB[s]\ \text{in } I_D/I_D^2,\qquad
\kappa'=\kappa-\frac{\log(AB)}{2\pi i}\quad\text{in the nonseparating case}.
\end{equation}
Here \(\kappa\) is the boundary value of \(\Omega_{22}-(2\pi i)^{-1}\log s\)
in the marking of Proposition~\ref{prop:two-period-limits}; the logarithms are
chosen compatibly.
\end{lemma}
\begin{proof}
Near the node, \(u,u'\) define the same reduced component of
\(\pi^{-1}D\), and \(v,v'\) define the other. Hence
\(u'=uU\), \(v'=vV\) for holomorphic units \(U,V\) whose values at the
node are \(A,B\). Since \(s'/s=UV\), its restriction to \(D\) is
\(AB\). The branch tangent bases are multiplied by \(A^{-1},B^{-1}\);
these factors cancel \(AB\) in \(\upsilon_\xi^{\mathrm{sm}}\).
Taking logarithms gives \eqref{eq:plumbing-coordinate-change}.
\end{proof}

In the nonseparating case let \(X=\widetilde C_0\). Once the elliptic differential
\(\alpha_1\) is chosen, take the branch coordinates to satisfy
\begin{equation}\label{eq:normalized-plumbing}
\dd u(\partial_O)=\dd v(\partial_P)=1,
\end{equation}
where \(\alpha_1(\partial_O)=\alpha_1(\partial_P)=1\) as in
Definition~\ref{def:normalized-tangents}. Two such choices have \(A=B=1\), so
\eqref{eq:plumbing-coordinate-change} gives the same \(\kappa\) for the fixed marking and logarithms.

Choose symplectic cycles $(A_1,A_2,B_1,B_2)$ and normalize the
holomorphic differentials by $\int_{A_i}\omega_j=\delta_{ij}$.
Write $\Omega_{ij}=\int_{B_i}\omega_j$.
\begin{proposition}\label{prop:two-period-limits}
The following expansions have holomorphic remainders.
\begin{enumerate}
\item Suppose the node is nonseparating, so that \(C_0=X/(O\sim P)\) with \(X\) elliptic. Orient \(A_2\) positively about \(P\). Outside the plumbing annulus, \(B_2\) runs from the boundary of the coordinate disk at \(O\) to that at \(P\); it returns through the annulus. Let \(\gamma\) be the induced path from \(O\) to \(P\) on the normalization \(X\), and denote the surviving cycles there by \(A_1,B_1\). Let \(\omega_X\) be the unique holomorphic differential on \(X\) satisfying \(\int_{A_1}\omega_X=1\). Put \(\tau=\int_{B_1}\omega_X\) and \(z=\int_\gamma\omega_X\). Then
\begin{equation}\label{eq:nonsep-period-matrix}
\Omega(s)=\begin{pmatrix}
\tau&z\\ z&\dfrac{\log s}{2\pi i}+\kappa
\end{pmatrix}+O(s).
\end{equation}
Let \(\eta_{P-O}\) be the unique meromorphic differential on \(X\) with simple poles only at \(O,P\), respective residues \(-1,1\), and zero
\(A_1\) period.  Compute \(I_{u,v}(\eta_{P-O})\) along \(\gamma\) by \eqref{eq:reg}, with endpoint logarithms compatible with the return of \(B_2\) through the annulus. Then \(2\pi i\kappa=I_{u,v}(\eta_{P-O})\).
\item Suppose the node separates two elliptic curves \(X_1,X_2\). Let \(Q_i\in X_i\)
be the preimage of the node. Keep the pointed curves \((X_i,Q_i)\) and their coordinates
\(u_i\) fixed as they are plumbed by \(u_1u_2=s\). Write
\(\omega_i=f_i(u_i)\dd u_i\) for the differential on \(X_i\) with \(A_i\) period one,
and put \(\tau_i=\int_{B_i}\omega_i\). Then
\begin{equation}\label{eq:separating-period-matrix}
\Omega(s)=\operatorname{diag}(\tau_1,\tau_2)
-2\pi i s f_1(0)f_2(0)
\begin{pmatrix}0&1\\1&0\end{pmatrix}+O(s^2).
\end{equation}
Thus \(\epsilon=\Omega_{12}\) vanishes along \(s=0\) and
\(\partial_s\epsilon|_{s=0}=-2\pi i f_1(0)f_2(0)\ne0\).
On the marked chart, \(\epsilon\) may therefore be used in place of
\(s\) as the local parameter for smoothing the node.
\end{enumerate}
\end{proposition}

\begin{proof}
In the nonseparating case the two normalized differentials extend as sections of the relative
dualizing sheaf. On the normalization their limits are \(\omega_X\) and
\((2\pi i)^{-1}\eta_{P-O}\), respectively \cite[Section~5.1.2]{HN}.
Apply Lemma~\ref{lem:crossing-finite-part} first to \(\omega_1(s)\): its \(A_2\) period
is zero, so \(\Omega_{21}\to\int_O^P\omega_X=z\). Symmetry gives the same limit for
\(\Omega_{12}\). Applying the lemma to \(2\pi i\,\omega_2(s)\), whose \(A_2\) period
is \(2\pi i\), gives \(2\pi i\kappa=I_{u,v}(\eta_{P-O})\).
The holomorphic expansions and the first order term for fixed separating components
are the plumbing formulas of \cite[Sections~5.1.1--5.1.2]{HN} after division of the
differentials there by \(2\pi i\); see also \cite{Yamada,Taniguchi89,Taniguchi91}.
A nonzero holomorphic differential on an elliptic curve has no zeros, so both \(f_i(0)\) are nonzero and
\(\partial_s\Omega_{12}(0)\ne0\).
\end{proof}

\subsection{An explicit family of enhanced curves of genus two}

We construct a family in the enhanced genus two space \(\mathsf T_2\) of \cite{CMY} whose nodal limit has normalization \(X\) and points \(O,P\) over the node. Work on the finite \'etale cover on which a point \(R\) with \(2R=P\) is chosen. Translation by \(-R\) sends \((O,P)\) to \((-R,R)\). Write \(R=(A,B)\), so that \(-R=(A,-B)\). Since \(2R=P\ne O\), we have \(B\ne0\). The duplication formulas express \((a,b)=(x(P),y(P))\) in these coordinates:
\begin{equation}\label{eq:halfpoint}
 a=\frac{(6A^2-t_2/2)^2}{4B^2}-2A,\qquad
 b=-B+\frac{6A^2-t_2/2}{B}(A-a).
\end{equation}

Split the double root at \(X=A\) by replacing \((X-A)^2\) with \((X-A)^2-s\). This gives
\begin{equation}\label{eq:explicit-quintic}
 Y^2=((X-A)^2-s)(4X^3-t_2X-t_3).
\end{equation}
Since \(\Delta\ne0\) and \(4A^3-t_2A-t_3=B^2\ne0\), the quintic has five simple
roots for sufficiently small \(s\ne0\); its smooth projective model \(C_s\) has genus two.
At \(s=0\), the normalization map is
\[
E:\ \widetilde y^{\,2}=4X^3-t_2X-t_3,\qquad
(X,\widetilde y)\longmapsto (X,(X-A)\widetilde y).
\]
It identifies \(-R=(A,-B)\) with \(R=(A,B)\).

Set \(u=X-2A/5\) and \(v=Y/2\). The equation becomes
$ v^2=u^5+p_2(s)u^3+p_3(s)u^2+p_4(s)u+p_5(s).$
On the smooth fibers, consider the holomorphic frame
\[
 \beta_1=(X-A)\frac{\dd X}{Y},
 \qquad
 \beta_2=B\frac{\dd X}{Y}.
\]
The limit of \(\beta_1\) is \(\dd X/\widetilde y\), while the limit of
\(\beta_2\) has residues \((-1,1)\) at \((-R,R)\). With respect to the
basis \((\dd u/v,u\dd u/v)\), and with rows indexed by
\(\beta_1,\beta_2\), the frame matrix and determinant are
\begin{equation}\label{eq:quintic-frame}
 S_1=\begin{pmatrix}
 -\frac{3A}{10}&\frac{1}{2}\\[2pt]
 \frac{B}{2}&0
 \end{pmatrix},
 \qquad
 \delta=\det S_1=-\frac B4.
\end{equation}
The chart of \cite{CMY} used below requires \(p_5(s)\ne0\). In fact \(p_5(s)=(9A^2/25-s)g(2A/5)\), where \(g(X)=X^3-t_2X/4-t_3/4\). Thus if \(g(2A/5)=0\), no fiber obtained by varying \(s\) alone lies in this chart. The frame itself is defined for \(B\Delta\ne0\); at points outside the chart the formulas for the generators are continued by varying the boundary parameters as well.

\begin{lemma}
\label{lem:quintic-normal-parameter}
The node of the central fiber of \eqref{eq:explicit-quintic} is \((X,Y)=(A,0)\).
Since \(B\ne0\), there is a unique holomorphic germ \(G\) at \(A\) with
\(G^2=4X^3-t_2X-t_3\) and \(G(A)=B\). The local coordinates
\[
 u_-=X-A-Y/G(X),\qquad u_+=X-A+Y/G(X)
\]
satisfy $u_-u_+=s$.  At the points $-R,R$ of the normalization, the coordinates
$\widetilde u_-=-u_-/(2B)$ and $\widetilde u_+=u_+/(2B)$
have tangent vectors dual to $\alpha_1$.  Their product is
\begin{equation}\label{eq:quintic-normal-parameter}
 \rho=-\frac{s}{4B^2}.
\end{equation}
\end{lemma}
\begin{proof}
The identity $(Y/G(X))^2=(X-A)^2-s$ gives $u_-u_+=s$.
At $-R$ and $R$, respectively,
$\dd u_-=2\dd X$, $\dd u_+=2\dd X$, and
$\alpha_1=-\dd X/B$, $\alpha_1=\dd X/B$.
The rescaled coordinates therefore have differential $\alpha_1$
at the corresponding point.
\end{proof}

Put $\omega_j=u^{j-1}\dd u/v$ for $1\le j\le4$ and set
\begin{equation}\label{eq:CMY-full-frame}
\gamma=S(s)\omega,\qquad
S(s)=\frac{1}{10B}
\begin{pmatrix}
 -3AB & 5B & 0 & 0\\
 5B^2 & 0 & 0 & 0\\
 0 & 0 & 5B & 0\\
 0 & 5p_2(s) & 3A & 15
\end{pmatrix}.
\end{equation}
\begin{proposition}\label{prop:CMY-frame}
The classes $\gamma$ extend across $s=0$ to a frame of the
logarithmic de Rham bundle, with $(\gamma_1,\gamma_2)=(\beta_1,\beta_2)$
and pairing matrix $J=\left(\begin{smallmatrix}0&I_2\\-I_2&0\end{smallmatrix}\right)$. For $s \ne 0$ with $p_5(s)\ne0$, the enhanced curve
$(C_s,\infty,\gamma)$ defines a map to $\mathsf T_2$.
\end{proposition}
\begin{proof}
In the basis $\omega$, the cup product matrix \cite[(3.10)]{CMY} is
\[
 \Psi(s)=\begin{pmatrix}
 0&0&0&\dfrac43\\[3pt]
 0&0&4&0\\[3pt]
 0&-4&0&\dfrac{4p_2(s)}3\\[3pt]
 -\dfrac43&0&-\dfrac{4p_2(s)}3&0
 \end{pmatrix}.
\]
Multiplication gives $S\Psi S^{\tr}=J$. The first two classes are
$\beta_1,\beta_2$, and their determinant $-B/4$ is nonzero; hence
they form the required holomorphic frame on every smooth fiber.

Near the node, Lemma~\ref{lem:quintic-normal-parameter} gives
\(\dd X/Y=G(X)^{-1}\dd u_+/u_+\) along the fibers, so each \(\omega_j\)
has a logarithmic representative. At infinity take \(u=t^{-2}\) and
\(v=t^{-5}H(t,s)\), where \(H(0,s)=1\) and
\(H^2=1+p_2(s)t^4+p_3(s)t^6+p_4(s)t^8+p_5(s)t^{10}\).
Then \(\omega_j=-2t^{4-2j}H^{-1}\dd t\). Since \(H\) is even in \(t\),
there is no residue term. Integrating the principal part gives a meromorphic
primitive depending holomorphically on the parameters. Subtracting its
differential near infinity gives a \v Cech representative of the extended class.
The pairing matrix is still \(J\) at \(s=0\), so the four extended classes form a frame.
\end{proof}

We now choose a basis whose limits can be compared with the elliptic graph, holomorphic, second kind and residue classes, in that order. Put \(r_0=9A^2/(25B)\), \(l_0=3A(84A^2-25t_2)/(500B^2)\), and set
\begin{equation}\label{eq:explicit-mixed-frame}
e_0=-\gamma_4+r_0\gamma_1+l_0\gamma_2,\quad
e_1=\gamma_1,\quad
e_2=\gamma_3-\tfrac A5\gamma_1-r_0\gamma_2,\quad e_3=\gamma_2 .
\end{equation}
\begin{proposition}\label{prop:explicit-graph-frame}
The pairing matrix of \((e_0,e_1,e_2,e_3)\) is \(\Phi\) of \eqref{eq:Phi}.
The limiting classes have the following meromorphic representatives on \(E\), in this order:
\begin{equation}\label{eq:normalization-images}
-\dd\!\left(\frac{\widetilde y}{2B}\right),\qquad
\frac{\dd X}{\widetilde y},\qquad
\frac{X\,\dd X}{\widetilde y},\qquad
\frac{B\,\dd X}{(X-A)\widetilde y}.
\end{equation}
The first differential represents \(\alpha_0\) in \(H^1_{\dR}(E,\{-R,R\})\).
\end{proposition}
\begin{proof}
The pairing follows by substituting
\eqref{eq:explicit-mixed-frame} into $J$.
At $s=0$, polynomial division gives
\[
 \gamma_3=\frac{X\dd X}{\widetilde y}+\frac A5\gamma_1+r_0\gamma_2,
 \qquad
 \gamma_4=\dd\!\left(\frac{\widetilde y}{2B}\right)
                          +r_0\gamma_1+l_0\gamma_2.
\]
Here $\dd\widetilde y=(6X^2-t_2/2)\dd X/\widetilde y$.
These identities give \eqref{eq:normalization-images}, using the primitives at infinity in the proof of Proposition~\ref{prop:CMY-frame}. The function \(-\widetilde y/(2B)\) takes the values \(1/2\) at \(-R\) and \(-1/2\) at \(R\). Its differential therefore has zero periods on closed cycles and integral \(-1\) from \(-R\) to \(R\). These periods characterize \(\alpha_0\) in relative cohomology.
\end{proof}


Write \(E=\C/\Lambda\) with coordinate \(w\) on its universal cover. Choose a lift
\(\widetilde z_R\) of \(R\). The translation \(z_E=w+\widetilde z_R\) sends the
lifts \(-\widetilde z_R,\widetilde z_R\) of the two branches to \(0,\widetilde z_P\),
where \(\widetilde z_P=2\widetilde z_R\) is the lift of \(P\) used in the periods below.
Put \(\nu=-(6A^2-t_2/2)/(2B)\).
\begin{lemma}\label{lem:exact-limiting-frame}
Under this translation the limiting de Rham classes are
\begin{equation}\label{eq:explicit-boundary-coordinates}
 (e_0,e_1,e_2,e_3)=
 (\alpha_0,\alpha_1,\beta+\nu\alpha_0,\eta_P+\nu\alpha_1).
\end{equation}
The last equality also holds for stable differentials on the
normalization.
\end{lemma}
\begin{proof}
On the normalization $X=\wp_\Lambda(w)$,
$\widetilde y=\wp_\Lambda'(w)$ and the points over the node have lifts
$-\widetilde z_R,\widetilde z_R$.  The zeta addition formula gives
\[
 \frac{\wp_\Lambda'(\widetilde z_R)}
      {\wp_\Lambda(w)-\wp_\Lambda(\widetilde z_R)}
 =\zeta_\Lambda(w-\widetilde z_R)-\zeta_\Lambda(w+\widetilde z_R)
                                      +2\zeta_\Lambda(\widetilde z_R).
\]
The duplication formula gives
$2\zeta_\Lambda(\widetilde z_R)-\zeta_\Lambda(\widetilde z_P)=\nu$.
After $z_E=w+\widetilde z_R$ these identities prove
$e_3=\eta_P+\nu\alpha_1$ as stable differentials, hence also
as classes in the logarithmic complex.

We have \(e_0=\alpha_0\) by Proposition~\ref{prop:explicit-graph-frame} and \(e_1=\dd w=\dd z_E=\alpha_1\). On each marked elliptic cycle, \(e_2=\wp_\Lambda(w)\dd w\) and \(\beta\) have the same period \(-\eta_j\); hence their difference is a multiple of \(\alpha_0\). The primitive \(-\zeta_\Lambda\) gives the relative period \(-2\zeta_\Lambda(\widetilde z_R)\) for \(e_2\) and the finite part \(-\zeta_\Lambda(\widetilde z_P)\) for \(\beta\) in the normalization of Section~\ref{sec:moduli}. By the zeta duplication formula their difference is \(-\nu\). As the relative period of \(\alpha_0\) is \(-1\), we obtain \(e_2-\beta=\nu\alpha_0\).
\end{proof}

We now introduce the relative frame parameters $c,t_1$. Set
$r=c-\nu$ and order the symplectic de Rham frame as follows:
\begin{equation}\label{eq:CMY-relative-parameters}
 \bigl(e_1,\ e_3+re_1,\ e_2+t_1e_1+re_0,\ -e_0\bigr).
\end{equation}
In the order \eqref{eq:CMY-relative-parameters} the pairing matrix is \(J\), and
the first two classes span \(H^0(C_s,\Omega^1_{C_s})\) for \(s\ne0\).
For the elliptic comparison, reorder the classes and replace \(-e_0\) by \(e_0\), obtaining
\((e_0,e_1,e_2+t_1e_1+re_0,e_3+re_1)\).
Substituting \eqref{eq:explicit-boundary-coordinates} and \(r=c-\nu\)
gives precisely \eqref{eq:frame} with \(d=0\), using the relative classes of
Section~\ref{sec:moduli} also when \(a=0\).

\subsection{The limiting period map and connection}
At a nonseparating node, the only row of the period matrix that can diverge is the row indexed by \(B_2\). In the plumbing \(u_Ou_P=\rho\), represent \(B_2\) by a path from the \(O\) end to the \(P\) end outside the annulus, followed by a return through the annulus. Orient the vanishing cycle \(A_2\) positively around \(P\). The following lemma identifies the finite part of the \(B_2\) period with an integral on the normalization \(X\). Use the induced path from \(O\) to \(P\) and logarithms compatible with the return through the annulus.
\begin{lemma}\label{lem:crossing-finite-part}
Let \(\xi\) be a holomorphic family of meromorphic differentials, regular as a section of the relative dualizing sheaf
near the node, with no other pole on the crossing path. Suppose that its restriction \(\xi(0)\) to the normalization has
residues \(-r,r\) at \(O,P\), respectively. Then
\[
L_\xi(\rho)=\int_{B_2}\xi(\rho)
-\frac{\log\rho}{2\pi i}\int_{A_2}\xi(\rho)
\]
extends holomorphically across \(\rho=0\). Use the endpoint tangents $\partial_O,\partial_P$ dual to $\alpha_1$,
so that $\dd u_O(\partial_O)=\dd u_P(\partial_P)=1$.
With these tangents and the logarithms fixed above, \eqref{eq:reg} gives
\begin{equation}\label{eq:crossing-finite-part}
L_\xi(0)=I_{u_O,u_P}(\xi(0))=\operatorname{Reg}\int_O^P\xi(0).
\end{equation}
\end{lemma}
\begin{proof}
Put \(u=u_O,\ v=u_P\). Near the node write \(\xi=f(u,v)\,\dd v/v\), with
\(f=\sum_{m,n\ge0}c_{mn}u^mv^n\); the coefficients are holomorphic in the remaining base parameters.  Define
\[
 a(\rho)=\sum_{n\ge0}c_{nn}\rho^n,\qquad
  H(u,v)=\sum_{\substack{m,n\ge0\\m\ne n}}
             \frac{c_{mn}}{n-m}u^mv^n.
\]
On a smooth fiber \(uv=\rho\), one has \(\dd u/u=-\dd v/v\). Thus, if \(\dd_{\mathrm{rel}}\) denotes differentiation
along the fiber,
\(\xi=a(\rho)\,\dd v/v+\dd_{\mathrm{rel}}H\).
In particular, \(\int_{A_2}\xi=2\pi i\,a(\rho)\) and \(a(0)=r\).
Fix small nonzero $\epsilon_O,\epsilon_P$. Let $\gamma_\rho$ be the part of $B_2$ outside the annulus, oriented from the boundary point $u=\epsilon_O$ near $O$ to the boundary point $v=\epsilon_P$ near $P$. The remaining part of \(B_2\) runs in the annulus from \((u,v)=(\rho/\epsilon_P,\epsilon_P)\) to \((\epsilon_O,\rho/\epsilon_O)\);
with the chosen logarithms its integral is
\[
 a(\rho)(\log\rho-\log\epsilon_O-\log\epsilon_P)
 +H\left(\epsilon_O,\frac{\rho}{\epsilon_O}\right)
 -H\left(\frac{\rho}{\epsilon_P},\epsilon_P\right).
\]
For fixed \(\epsilon_O,\epsilon_P\), the displayed annulus integral minus \(a(\rho)\log\rho\) is holomorphic at \(\rho=0\). So is the integral along \(\gamma_\rho\): its endpoints are holomorphic sections, and along this path away from the poles, local primitives can be chosen holomorphically in the base parameters. Their sum is \(L_\xi(\rho)\), independently of the cutoffs. Set \(\rho=0\) before letting the cutoffs tend to zero. Then \(H(\epsilon_O,0)-H(0,\epsilon_P)\to0\), and the limit of \(\int_{\gamma_0}\xi(0)-r\log\epsilon_O-r\log\epsilon_P\) is the regularized integral in \eqref{eq:reg}.
\end{proof}

Return to the family \eqref{eq:explicit-quintic}. The five coordinates on its boundary \(\mathcal B\) are
\((A,c,t_1,t_2,t_3)\), with a local choice of \(B\) satisfying \(B^2=4A^3-t_2A-t_3\) and \(B\Delta\ne0\). Doubling
\(R=(A,B)\) gives \(P=(a,b)\) by \eqref{eq:halfpoint}. Thus
\(\ell_0:\mathcal B\to\T_0^{\mathrm{an}}\) sends \((A,c,t_1,t_2,t_3)\) to \((a,b,c,t_1,t_2)\). The smoothing parameter \(\rho=-s/(4B^2)\) gives the sixth coordinate on \(\mathcal B\times\mathbb D\).

Form the period matrix \(P_2\) with rows indexed by \((B_2,B_1,A_1,A_2)\), using on the normalization the lifts of the
elliptic cycles and the path fixed in \eqref{eq:Q}. Its columns correspond to the frame
\[
\bigl(e_0,\ e_1,\ e_2+t_1e_1+(c-\nu)e_0,\ e_3+(c-\nu)e_1+c_{\mathrm{per}}e_0\bigr).
\]
The last column includes the same constant shift \(c_{\mathrm{per}}e_0\) as in
Proposition~\ref{prop:normalized-period-matrix}.
To compute this column, first take the vanishing and crossing periods of \(e_3\) alone:
\[
a_v(\rho)=\frac1{2\pi i}\int_{A_2} e_3,\qquad
\mathcal R=\lim_{\rho\to0}\left(\int_{B_2}e_3-a_v(\rho)\log\rho\right).
\]
For the same path and logarithms, let \(\theta\) be the coefficient for which the regularized relative period of \(\eta_P+\nu\alpha_1+\theta\alpha_0\) equals \(\mathcal R\). Since the relative period of \(\alpha_0\) is \(-1\), equation~\eqref{eq:relative-third} gives
\begin{equation}\label{eq:central-limit-explicit}
 \theta=-2\log\sigma_\Lambda(\widetilde z_P)
       +\bigl(\nu+\zeta_\Lambda(\widetilde z_P)\bigr)\widetilde z_P
       -\chi-\mathcal R,
\end{equation}
where $\chi=\log\sigma_\Lambda(-\widetilde z_P)
-\log\sigma_\Lambda(\widetilde z_P)\in(2\Z+1)\pi i$
is the branch term of \eqref{eq:relative-third}.

Let \(N\) be the logarithm of positive monodromy on the four cycle rows, numbered \(0,1,2,3\), and fix the branch of \(\log\rho\) used in the crossing integral. Remove the logarithmic term by setting
\begin{equation}\label{eq:regularized-genus-two-period}
\widehat P_2=\exp\!\left(-\frac{\log\rho}{2\pi i}N\right)P_2.
\end{equation}
\begin{theorem}\label{thm:enhanced-nonsep-limit}
For the family and frame above, the following hold.
\begin{enumerate}
\item $N=E_{03}$. The function $a_v$ and the matrix $\widehat P_2$ extend holomorphically
to $\rho=0$. Write $v=\widetilde z_P$, $g=c+\zeta_\Lambda(v)$ and
$w_0=-2\log\sigma_\Lambda(v)+gv-2\log(2\pi k)$, with the branches
of \eqref{eq:Q}. The limiting matrix is invertible and is given by
\begin{equation}\label{eq:enhanced-period-limit}
\widehat P_2(0)=Q(a,b,c,t_1,t_2,0)=
\begin{pmatrix}
 -1&v&t_1v-g&w_0\\
 0&\omega_1&t_1\omega_1-\eta_1&g\omega_1-\eta_1v\\
 0&\omega_2&t_1\omega_2-\eta_2&g\omega_2-\eta_2v\\
 0&0&0&2\pi i
\end{pmatrix}.
\end{equation}
Moreover,
\begin{equation}\label{eq:explicit-central-zero}
 a_v(0)=1,\qquad \theta=0 .
\end{equation}
\item With \(D_0=\operatorname{diag}(1,k,k,k^2)\) as in \eqref{eq:Q}, the map
\(\ell=\operatorname{Per}^{-1}([\widehat P_2(0)D_0]):\mathcal B\to\T^{\times,\mathrm{an}}\)
lifts \(\ell_0\) and satisfies \(\ell^*h=1\).
\end{enumerate}
\end{theorem}

\begin{proof}
Positive monodromy sends \(B_2\) to \(B_2+A_2\) and fixes the other cycles, so \(N=E_{03}\). Consequently \eqref{eq:regularized-genus-two-period} subtracts \((\log\rho)/(2\pi i)\) times the last row of \(P_2\) from its first row. Each column satisfies the hypotheses of Lemma~\ref{lem:crossing-finite-part} by Proposition~\ref{prop:CMY-frame}. The lemma extends the first row holomorphically; its residue calculation extends the last. The two elliptic rows extend
by integration along cycles away from the node.

It remains to identify these limiting periods. By Lemma~\ref{lem:exact-limiting-frame}, \(e_3\) restricts to \(\eta_P+\nu\alpha_1\), with residue \(1\) at \(P\). Thus \(a_v(0)=1\); by \eqref{eq:crossing-finite-part}, \(\mathcal R\) is the regularized integral of that differential, so \(\theta=0\). Equation \eqref{eq:explicit-boundary-coordinates} identifies the first three columns with the periods of \(\alpha_0,\alpha_1,\alpha_2\). Adding \((c-\nu)e_1+c_{\mathrm{per}}e_0\) to \(e_3\) gives the limiting fourth class \(\eta_P+c\alpha_1+c_{\mathrm{per}}\alpha_0\); its finite crossing period is \(\mathcal R+(c-\nu)\widetilde z_P-c_{\mathrm{per}}\). These are exactly the four columns of \(Q(a,b,c,t_1,t_2,0)\) in \eqref{eq:Q}. The pairing identity \eqref{eq:period-cup-identity} gives invertibility. Theorem~\ref{thm:inverse} identifies this normalized boundary period with the elliptic enhancement having \(d=0\), hence \(h=1\).
\end{proof}
For the holomorphic differentials on a smooth fiber we continue to use \(e_1,e_3+(c-\nu)e_1\) from \eqref{eq:CMY-relative-parameters}. The term \(c_{\mathrm{per}}e_0\) enters only the period comparison.

\begin{theorem}\label{thm:connection-degeneration}
Fix the parameter \(\rho=-s/(4B^2)\), put \(Q=\widehat P_2|_{\mathcal B}\), and use the frame of
Theorem~\ref{thm:enhanced-nonsep-limit}.
\begin{enumerate}
\item The Gauss-Manin connection has the logarithmic matrix
\begin{equation}\label{eq:connection-decomposition}
\bigl(P_2^{-1}\dd P_2\bigr)^{\tr}
=\bigl(\widehat P_2^{-1}N\widehat P_2\bigr)^{\tr}
\frac{\dd\rho}{2\pi i\rho}
+\bigl(\widehat P_2^{-1}\dd\widehat P_2\bigr)^{\tr}.
\end{equation}
Its residue $(Q^{-1}NQ)^{\tr}/(2\pi i)=-E_{30}$ is nilpotent of rank one.
\item Pull back the holomorphic second term of  \eqref{eq:connection-decomposition} by the inclusion \(\mathcal
B=\{\rho=0\}\hookrightarrow\mathcal B\times\mathbb D\). With \(\dd_{\mathcal B}\) denoting differentiation on the
boundary, the resulting matrix is
\begin{equation}\label{eq:boundary-connection-A}
\bigl(Q^{-1}\dd_{\mathcal B}Q\bigr)^{\tr}=\ell^*\A .
\end{equation}
Replacing $Q$ by $QD_0$ gives the boundary matrix $D_0(\ell^*\A)D_0^{-1}$.
\end{enumerate}
\end{theorem}
\begin{proof}
Write \(P_2=L\widehat P_2\), where \(L=\exp((\log\rho)N/(2\pi i))\).
The product rule and \(L^{-1}\dd L=N\dd\rho/(2\pi i\rho)\) give
\eqref{eq:connection-decomposition}. Since \(\widehat P_2\) is holomorphic
and invertible along the boundary, the residue is \((Q^{-1}NQ)^{\tr}/(2\pi i)\).
Let \(b_0,\ldots,b_3\) be the standard basis columns. By \eqref{eq:Q},
\(Qb_0=-b_0\) and \(b_3^{\tr}Q=2\pi i b_3^{\tr}\). Thus
\(Q^{-1}NQ=Q^{-1}b_0b_3^{\tr}Q=-2\pi i E_{03}\), giving residue \(-E_{30}\).
Theorem~\ref{thm:enhanced-nonsep-limit} identifies \(Q\) with the elliptic period matrix of \(\ell\). Differentiating on
\(\mathcal B\) and applying Proposition~\ref{prop:analytic-connection} gives \eqref{eq:boundary-connection-A}. In
particular, \(\ell^*(\dd h/h)=0\) since \(\ell^*h=1\). The change \(Q\mapsto QD_0\) gives the stated conjugation because
\(D_0\) is constant and diagonal.
\end{proof}

\subsection{Restrictions of genus two generators}
We now evaluate the generators of \cite{CMY} on the enhanced curve \((C_s,\infty,\gamma)\), using the frame
\(\gamma=S(s)\omega\) of \eqref{eq:CMY-full-frame}. Thus the frame entries in the formulas below are those of \(S(s)\),
before the changes of basis used to form \(P_2\).
The coefficients of the quintic depend on the chosen half \(R\) of \(P\). We first compute the generator values for that
choice and then identify combinations independent of it.

At $s=0$, write $p_i=p_i(0)$ for the coefficient of $u^{5-i}$
($2\le i\le5$) in the right hand side of \eqref{eq:explicit-quintic}
after setting $X=u+2A/5$ and dividing by $4$. Expansion gives
\begin{equation}\label{eq:quintic-coefficients}
\begin{aligned}
p_2&=-\frac{3A^2}{5}-\frac{t_2}{4}, &
p_3&=-\frac{2A^3}{25}+\frac{At_2}{5}-\frac{t_3}{4},\\
p_4&=\frac{12A^4}{125}+\frac{3A^2t_2}{100}+\frac{3At_3}{10}, &
p_5&=\frac{72A^5}{3125}-\frac{9A^3t_2}{250}-\frac{9A^2t_3}{100}.
\end{aligned}
\end{equation}
\begin{proposition}\label{prop:CMY-explicit}
On the chart \(p_5(s)\ne0\) over \(S\setminus D\), the pullback of each generator in \cite[(4.11)]{CMY} is a polynomial in the \(p_i(s)\) and the entries of \(S(s)\), divided by a power of \(\delta=-B/4\). It therefore extends holomorphically to \(s=0\) when \(B\Delta\ne0\), with boundary value obtained by substituting \(p_i=p_i(0)\) and \(S=S(0)\). In particular,
\begin{equation}\label{eq:CMY-five-limits}
\begin{aligned}
T_4&=-4p_2/B, &T_8&=16p_4/B^2, &T_{12}&=-64p_3^2/B^3,\\
T_{16}&=256p_3p_5/B^4, &T_{20}&=-1024p_5^2/B^5.
\end{aligned}
\end{equation}
In the full frame \eqref{eq:CMY-full-frame}, the additional boundary values are
\[
U_{3,11}=4p_3/B^2,\qquad U_{5,11}=-16p_5/B^3,\qquad
U_{j,44}=(p_2/B)^2 U_{j,11}\quad(j=3,5).
\]
\end{proposition}
\begin{proof}
By \cite[(4.11)]{CMY}, the generators \(T_4,T_8,T_{12},T_{16},T_{20}\) are, respectively,
\[
\frac{p_2(s)}{\delta},\quad \frac{p_4(s)}{\delta^2},\quad
\frac{p_3(s)^2}{\delta^3},\quad \frac{p_3(s)p_5(s)}{\delta^4},\quad
\frac{p_5(s)^2}{\delta^5}.
\]
For \(i=1,4\), the remaining generators in the statement are \(U_{3,ii}=p_3(s)s_{i2}^2/\delta^2\) and
\(U_{5,ii}=p_5(s)s_{i2}^2/\delta^3\), with \(s_{ij}=S(s)_{ij}\). All denominators are nonzero because \(\delta=-B/4\).
The formulas therefore extend from the chart in the full family, even if a fixed one parameter slice does not meet the
chart. Setting \(s=0\), \(s_{12}=1/2\), and \(s_{42}(0)=p_2/(2B)\) gives the asserted values.
\end{proof}

Define the following polynomials in the genus two generators:
\begin{align*}
 \mathcal A&=-3T_4^2-20T_8,\qquad
 \mathcal B_{\mathrm{inv}}=-3T_4T_{12}+9T_4^2T_8-20T_8^2+75T_{16},\\
 \mathcal C
 &= 12T_4^3T_{12}+18T_{12}^2-36T_4^4T_8
    -13T_4T_{12}T_8\\
 &\quad-88T_4^2T_8^2+160T_8^3-165T_4^2T_{16}
          -800T_8T_{16}-1125T_4T_{20}.
\end{align*}
Let \(\mathcal D\) be the discriminant in \(u\) of
\(u^5+p_2(s)u^3+p_3(s)u^2+p_4(s)u+p_5(s)\), divided by \(5^5\delta^{10}\).
\begin{corollary}\label{cor:invariant-boundary-values}
The following three combinations restrict to functions of the marked elliptic curve and its chosen differential, independently of the choice of \(R\) with \(2R=P\):
\begin{equation}\label{eq:CMY-invariant-limits}
\mathcal A|_{s=0}=-48a,\qquad
\mathcal B_{\mathrm{inv}}|_{s=0}=12t_2,\qquad
(4\mathcal A\mathcal B_{\mathrm{inv}}-3\mathcal C)|_{s=0}=864t_3.
\end{equation}
At $s=0$, the normalized discriminant $\mathcal D$ and its product with $\mathcal A$ have the expansions
\begin{equation}\label{eq:CMY-discriminant-leading}
\mathcal D=-\frac{1024\Delta}{3125B^2}s+O(s^2),\qquad
\mathcal A\mathcal D=\frac{49152a\Delta}{3125B^2}s+O(s^2).
\end{equation}
\end{corollary}
\begin{proof}
Substituting \eqref{eq:CMY-five-limits} into \(\mathcal A,\mathcal B_{\mathrm{inv}},4\mathcal A\mathcal B_{\mathrm{inv}}-3\mathcal C\) and using \(B^2=4A^3-t_2A-t_3\) and \eqref{eq:halfpoint} gives \eqref{eq:CMY-invariant-limits}. Its right hand sides depend only on \((X,O,P,\alpha_1)\).

For the discriminant, write $g(X)=X^3-t_2X/4-t_3/4$. The discriminant of $((X-A)^2-s)g(X)$ is $4s\,\operatorname{disc}(g)\operatorname{Res}((X-A)^2-s,g)^2$. Here $\operatorname{disc}(g)=-\Delta/16$, and at $s=0$ the resultant is $g(A)^2=B^4/16$. Translation does not change the discriminant. Dividing by \(5^5(-B/4)^{10}\) gives the coefficient of \(s\) in \(\mathcal D\). Since \(\mathcal D(0)=0\), multiplication by \(\mathcal A(0)=-48a\) gives that of \(\mathcal A\mathcal D\).
\end{proof}


Since \(B\Delta\ne0\), equation~\eqref{eq:CMY-discriminant-leading} shows that \(\mathcal D^{-1}\) is regular on \(S\setminus D\) and has a simple pole along \(s=0\).

By \cite[Section~7]{CMY}, the functions \(\mathcal B_{\mathrm{inv}}\), \(4\mathcal A\mathcal B_{\mathrm{inv}}-3\mathcal C\), \(\mathcal D\) and \(\mathcal A\mathcal D\) are algebraic representatives of the Siegel modular forms \(E_4^{(2)}, E_6^{(2)}, \chi_{10}, \chi_{12}\), respectively, up to a nonzero constant factor in each case.

\subsection{The smoothing parameter and the effective tangent coordinate}
We compute the first derivative of the Siegel parameter \(q_2\) at the boundary.
We first use \(\rho=-s/(4B^2)\) and set \(h_\rho=\ell^*h=1\), as in
Theorem~\ref{thm:enhanced-nonsep-limit}. For the limiting differential \(e_1\),
write \(\omega_2,\omega_1\) for its periods over \(A_1,B_1\), respectively, and
\(V=\widetilde z_P\) for its integral along the chosen path from \(O\) to \(P\).
Put \(\tau=\omega_1/\omega_2\), \(z=V/\omega_2\). The limiting \(A_1\) period
of \(e_3\) is \(p=-\eta_2V+(\nu+\zeta_\Lambda(V))\omega_2\).
\begin{proposition}\label{prop:normal-unit-comparison}
For the Siegel parameter \(q_2=e^{2\pi i\Omega_{22}}\), the quotient
\(q_2/\rho\) extends to a holomorphic unit with boundary value
\begin{equation}\label{eq:normal-unit-explicit}
U_\rho:=\left.\frac{q_2}{\rho}\right|_{\rho=0}
=-\frac{4\pi^2h_\rho}{\omega_2^2F(\tau,z)^2}.
\end{equation}
Equivalently, in the algebraic parameter \(s\),
\begin{equation}\label{eq:algebraic-period-parameter}
q_2=\frac{\pi^2}{B^2\omega_2^2F(\tau,z)^2}\,s+O(s^2).
\end{equation}
\end{proposition}
\begin{proof}
To compute \(\Omega_{22}\), write \(\varpi_2=\ell_1 e_1+\ell_3e_3\) for the
holomorphic differential with \(A_1,A_2\) periods \(0,1\). The matrix of \(A\) periods of \((e_1,e_3)\) extends holomorphically with invertible limit \(\left(\begin{smallmatrix}\omega_2&p\\0&2\pi i\end{smallmatrix}\right)\). Thus \(\ell_1,\ell_3\) are holomorphic up to the boundary and \((\ell_1(0),\ell_3(0))=(-p/\omega_2,1)/(2\pi i)\). For \(j=1,3\), let \(L_j\) be the finite part of the \(B_2\) period of \(e_j\) after subtracting \((2\pi i)^{-1}\log\rho\int_{A_2}e_j\). Lemma~\ref{lem:crossing-finite-part} gives \(L_1(0)=V\) and \(L_3(0)=\mathcal R\). Since \(\int_{A_2}\varpi_2=1\), the logarithmic term in \(2\pi i\Omega_{22}\) is exactly \(\log\rho\). Hence \(q_2/\rho=\exp(2\pi i(\ell_1L_1+\ell_3L_3))\) is a holomorphic unit. Restricting this unit to \(\rho=0\) and using
\eqref{eq:central-limit-explicit} and the value of \(p\), we obtain
\[
U_\rho=\left.\frac{q_2}{\rho}\right|_{\rho=0}
=\exp(\mathcal R-pV/\omega_2)
=-\frac{h_\rho}{\sigma_\Lambda(V)^2}
\exp\!\left(\frac{\eta_2V^2}{\omega_2}\right),
\]
where \(e^{-\chi}=-1\) and \(e^{-\theta}=h_\rho\). Using the expression for the Weierstrass sigma function in terms of the
odd Jacobi theta function \cite[Eq.~23.6.9]{DLMF} and
\(\partial_z\theta_1(\tau,0)=2\pi\eta(\tau)^3\), we obtain
\[
 \sigma_\Lambda(V)=\frac{\omega_2}{2\pi}
 \exp\!\left(\frac{\eta_2V^2}{2\omega_2}\right)F(\tau,z).
\]
This proves \eqref{eq:normal-unit-explicit}; substituting the fixed value of \(\rho\) gives \eqref{eq:algebraic-period-parameter}.
\end{proof}

We now allow the local defining function \(\rho\) of the divisor
\(D=\mathcal B\times\{0\}\subset S=\mathcal B\times\mathbb D\) to vary.
Let \(\pi:L^\times\to D\) be its normal bundle with the zero section removed.
Keep the branch order, homology marking, and de Rham frame of
Theorem~\ref{thm:enhanced-nonsep-limit} fixed. For each such \(\rho\), define
\[
Q_\rho=\left.\left(\exp\!\left(-\frac{\log\rho}{2\pi i}N\right)P_2\right)\right|_D,
\qquad \A_\rho=(Q_\rho^{-1}\dd_D Q_\rho)^{\tr}.
\]
The exponential removes the logarithmic term before restriction to \(D\); \(\dd_D\) denotes differentiation on \(D\). Let \(h_\rho\) be the \(h\) coordinate of \(\operatorname{Per}^{-1}([Q_\rho D_0])\), where \(D_0=\operatorname{diag}(1,k,k,k^2)\) and \(k=(-2\pi i)^{-1/2}\) is the fixed choice. A nonzero normal vector \(\xi\) over \(x\in D\) has coordinate \(r_\rho(\xi)=\dd\rho_x(\xi)\in\C^\times\).

\begin{theorem}
\label{thm:normal-bundle}\label{lem:normal-unit-choices}
For \(x\in D\) and \(0\ne\xi\in N_{D/S,x}\), let \(\ell_N(\xi)\)
denote the enhancement of \(\ell_0(x)\) with effective tensor
\(k^2\upsilon^{\mathrm{sm}}_\xi\). 
The regularized periods induce a holomorphic \(\C^\times\) equivariant
isomorphism over \(D\):
\[
 L^\times\xrightarrow{\sim}D\times_{\T_0^{\mathrm{an}}}\T^{\times,\mathrm{an}},
 \qquad \xi\longmapsto\bigl(\pi(\xi),\ell_N(\xi)\bigr).
\]
\begin{itemize}
\item
On \(L^\times\), with functions on \(D\) pulled back by \(\pi\), we have
\begin{equation}\label{eq:normal-line-comparison}
h_\xi:=(\ell_N^*h)(\xi)=r_\rho(\xi)h_\rho,\qquad
\ell_N^*\A=\pi^*\A_\rho-E_{30}\frac{\dd r_\rho}{r_\rho}.
\end{equation}
\item
Since \(q_2\) vanishes on \(D\), its differential at \(x\in D\)
descends to \(N_{D/S,x}=T_xS/T_xD\). With the chosen marking, its value is
\begin{equation}\label{eq:normal-line-q}
\dd q_2(\xi)=-\frac{4\pi^2 h_\xi}{\omega_2^2F(\tau,z)^2}.
\end{equation}
\end{itemize}
\end{theorem}

\begin{proof}
Write \(\rho'=v\rho\), where \(v\) is a holomorphic unit near \(D\), and put \(u=v|_D\). Then \(r_{\rho'}=ur_\rho\) and \(Q_{\rho'}=\exp(-(\log u)N/(2\pi i))Q_\rho\). For compatible logarithms these two changes cancel in
\[
Q(\xi)=\exp\!\left(\frac{\log r_\rho(\xi)}{2\pi i}N\right)Q_\rho.
\]
Changing a logarithm multiplies this matrix on the left by an integral monodromy transformation. Since \(N=E_{03}\) and the last row of \(Q_\rho\) is \((0,0,0,2\pi i)\), the exponential adds \(\log r_\rho(\xi)\) to the upper right entry of \(Q_\rho\) and fixes every other entry. Thus \(Q(\xi)D_0\in\Pi\), and \(\ell_N(\xi):=\operatorname{Per}^{-1}([Q(\xi)D_0])\) is independent of the choices of parameter and logarithms. By \eqref{eq:Q}, this change fixes the relative elliptic frame and gives \(h_\xi=r_\rho(\xi)h_\rho\). Multiplication by the nonzero function \(h_\rho\) has holomorphic inverse \(h\mapsto h/h_\rho\) in these fiber coordinates. This proves the equivariant isomorphism. In the normalized plumbing coordinates, \(h_\rho=1\) and the two branch tangents are \(\partial_O,\partial_P\); hence \eqref{eq:tensor} and Lemma~\ref{lem:plumbing-tensor} give the stated tensor. Differentiating \(Q(\xi)\) contributes the term \((Q_\rho^{-1}NQ_\rho)^{\tr}\dd\log r_\rho/(2\pi i)=-E_{30}\dd\log r_\rho\) in addition to \(\pi^*\A_\rho\), by Theorem~\ref{thm:connection-degeneration}. Finally \(\dd q_2(\xi)=U_\rho r_\rho(\xi)\). Under \(\rho'=v\rho\), both \(U_\rho\) and \(h_\rho\) are divided by \(u\). Their ratio is therefore given by \eqref{eq:normal-unit-explicit} for every \(\rho\), proving \eqref{eq:normal-line-q}.
\end{proof}

In particular, every regular function of index \(m\in\Z\) has the form
\(ph^m\) with \(p\in A_0\). On the punctured normal line at \(x\in D\) its pullback is \(\xi\mapsto p(\ell_0(x))h_\xi^m\), homogeneous of degree \(m\).

\begin{remark}
For a fixed family and marking, the limiting matrix depends only on the restriction of the de Rham frame to \(D\). Indeed, under \(P_2\mapsto MP_2G\), with \(M\) a constant change of marking and \(G\) a holomorphic frame change invertible along \(D\), one has \(N\mapsto MNM^{-1}\) and \(Q(\xi)\mapsto MQ(\xi)G_D\), where \(G_D=G|_D\). The resulting connection matrix is \(G_D^{\tr}(\ell_N^*\A)(G_D^{\tr})^{-1}+(G_D^{-1}\dd G_D)^{\tr}\), with \(G_D\) pulled back to \(L^\times\). In particular, in the same marking, two extensions that agree on \(D\) give the same limiting matrix. The parameter \(q_2\) is computed from the differentials with identity \(A\) period matrix, so replacing \(e_3\) by \(e_3+v e_1\) leaves \(q_2\) unchanged.
\end{remark}

\begin{corollary}\label{cor:normal-coefficient-conversion}
Fix the boundary coordinates and return to \(\rho=-s/(4B^2)\), so \(h_\rho=1\).
If \(f\) is meromorphic with poles supported on \(s=0\) and
\(f=s^m f_m+O(s^{m+1})\), with \(m\in\Z\) and \(f_m\) holomorphic
on \(\mathcal B\), then
\begin{equation}\label{eq:normal-coefficient-conversion}
f=q_2^m
\left(\frac{B^2\omega_2^2F(\tau,z)^2}{\pi^2}\right)^m f_m
+O(q_2^{m+1}).
\end{equation}
In particular, for the normalized quintic discriminant \(\mathcal D\) and its
product with \(\mathcal A=-3T_4^2-20T_8\) from
Corollary~\ref{cor:invariant-boundary-values}, we obtain
\begin{equation}\label{eq:discriminant-period-leading}
\bigl(\mathcal D,\,\mathcal A\mathcal D\bigr)
=\frac{1024\omega_2^2\Delta F(\tau,z)^2}{3125\pi^2}
(-1,48a)\,q_2+O(q_2^2).
\end{equation}
\end{corollary}
\begin{proof}
Invert \eqref{eq:algebraic-period-parameter}; its linear coefficient is nonzero
since \(B\omega_2F(\tau,z)\ne0\). Substitution gives
\eqref{eq:normal-coefficient-conversion}, and
\eqref{eq:CMY-discriminant-leading} gives \eqref{eq:discriminant-period-leading}.
\end{proof}


Let \(\mathfrak F\) be a scalar Siegel modular form of genus two and weight \(r\). To express \(\mathfrak F(\Omega(s))\) in the chosen algebraic basis, we compare the two bases of holomorphic differentials. Let \(f(s)\) be its algebraic representative in the chosen basis and write \(\mathcal P_A(s)\) for the matrix of \(A\) periods of that basis. Up to a fixed normalization constant \(c_f\),
\begin{equation}\label{eq:algebraic-classical-frame-conversion}
 \mathfrak F(\Omega(s))=c_f\det(\mathcal P_A(s))^r f(s).
\end{equation}
Indeed, the wedge of the chosen basis is \(\det\mathcal P_A\) times the wedge of the basis with identity \(A\) period matrix. For \((e_1,e_3)\) the limiting determinant is \(2\pi i\omega_2\). The replacement \(e_3\mapsto e_3+(c-\nu)e_1\) preserves this wedge; the additional term \(c_{\mathrm{per}}e_0\) appears only in the comparison of de Rham frames.

Corollary~\ref{cor:normal-coefficient-conversion} determines the leading coefficient with the coordinates on \(\mathcal B\) held fixed. Higher coefficients along the family also involve the variation of \(\Omega_{11}(s)\) and \(\Omega_{12}(s)\).

\subsection{Fourier--Jacobi coefficients}
Let \(\mathfrak F\) be a scalar holomorphic Siegel modular form of integral weight \(r\) for \(\Sp_4(\Z)\), with trivial character. Write
\(\Omega=\left(\begin{smallmatrix}\tau&z\\z&\omega\end{smallmatrix}\right)\)
and \(q_2=e^{2\pi i\omega}\). Its classical Fourier-Jacobi expansion \cite[Section~8]{vdG} converges normally on compact subsets of the Siegel upper half space:
\begin{equation}\label{eq:classical-FJ-new}
\mathfrak F(\Omega)=\sum_{m\ge0}\phi_m(\tau,z)q_2^m,\qquad
\phi_m\in J_{r,m}.
\end{equation}

\begin{corollary}\label{cor:GMCD-FJ}
\begin{enumerate}
\item[(i)]
On the classical period locus \((\tau,z)\in\mathcal H\),
each Fourier--Jacobi coefficient is the pullback of a regular function of index \(m\):
\begin{equation}\label{eq:GMCD-FJ-new}
\phi_m=\mathsf t^*(h^mp_m),\qquad
\mathfrak F(\Omega)=\sum_{m\ge0}(\mathsf t^*p_m)\bigl(F(\tau,z)^2 q_2\bigr)^m.
\end{equation}
Here \(p_m=p_m(a,b,t_2,t_3)\) has weight \(r+2m\) and is unique modulo
\(b^2=4a^3-t_2a-t_3\); its normal form is given in Corollary~\ref{thm:intro-FJ}.
The elliptic function \(\mathsf t^*p_m\) has poles only at lattice points,
of order at most \(2m\).
\item[(ii)] On the marked smoothing, put \(f(s)=\mathfrak F(\Omega(s))\).
If \(\phi_j\equiv0\) for \(0\le j<m\), the coefficient of order \(m\)
evaluated on a nonzero normal vector \(\xi\) is
\begin{equation}\label{eq:FJ-normal-leading}
\left.\frac{f}{\rho^m}\right|_D\bigl(\dd\rho(\xi)\bigr)^m
=\phi_m(\tau,z)\bigl(\dd q_2(\xi)\bigr)^m
=(k\omega_2)^r(-2\pi i)^m p_m(a,b,t_2,t_3)\,h_\xi^m .
\end{equation}
In the last expression, \((a,b,t_2,t_3)\) are the coordinates of the limiting elliptic curve with its chosen differential, evaluated at \(\pi(\xi)\).

\end{enumerate}
\end{corollary}
\begin{proof}
For \textup{(i)}, apply Proposition~\ref{prop:weak-ring} to each
\(\phi_m\) and use \(\mathsf t^*h=F^2\). The normal form gives uniqueness;
its monomials in \(a,b\) have pole order at most \(2m\).

For \textup{(ii)}, divide the Fourier-Jacobi expansion by \(q_2^m\).
The lower coefficients vanish identically, so the quotient extends
holomorphically to \(D\) with value \(\phi_m(\tau,z)\), where
\((\tau,z)=(\Omega_{11},\Omega_{12})|_D\). Since
\(\dd q_2(\xi)=(q_2/\rho)|_D\,\dd\rho(\xi)\), this proves the first
equality of \eqref{eq:FJ-normal-leading}, even at a zero of \(\phi_m\).

To compare the two elliptic frames, note that their period lattices are \(\omega_2(\Z\tau+\Z)\) for the limiting differential and \(k^{-1}(\Z\tau+\Z)\) on the classical locus. Hence a coordinate \(x\in\{a,b,t_2,t_3\}\) satisfies \(\mathsf t^*x=(k\omega_2)^{\wt(x)}x\), where the right hand side uses the limiting differential. Part~\textup{(i)} gives
\(\phi_m=F^{2m}(k\omega_2)^{r+2m}p_m\).
Substitute \eqref{eq:normal-line-q} and use \(-4\pi^2k^2=-2\pi i\).
\end{proof}

Put $\Delta_{\mathrm{ell}}=(E_4^3-E_6^2)/1728$.
Normalize the genus two Eisenstein series to have constant term $1$,
and choose scalar multiples $\widehat\chi_{10},\widehat\chi_{12}$
of the Igusa cusp forms so that their first Fourier-Jacobi
coefficients are $\Delta_{\mathrm{ell}}F^2$ and
$\Delta_{\mathrm{ell}}k^2\wp F^2$, respectively.

\begin{proposition}\label{prop:classical-leading}
Evaluate the elliptic coordinates along \(\mathsf t\), and hold
\((\tau,z)\in\mathcal H\) fixed as \(q_2\to0\).
\begin{enumerate}
\item[(1)]
The first terms of the four Fourier-Jacobi expansions are
\begin{equation}\label{eq:Igusa-leading}
\begin{aligned}
E_4^{(2)}&=12k^4t_2+O(q_2),&
E_6^{(2)}&=-216k^6t_3+O(q_2),\\
\widehat\chi_{10}&=-k^{12}\Delta h\,q_2+O(q_2^2),&
\widehat\chi_{12}&=-k^{12}\Delta ah\,q_2+O(q_2^2).
\end{aligned}
\end{equation}
Consequently, \(\widehat\chi_{12}/\widehat\chi_{10}\to a\).
\item[(2)]
The first coefficient \(\phi=-k^{12}\Delta h\)
of \(\widehat\chi_{10}\) and the coordinate \(a\) satisfy
\begin{equation}\label{eq:differential-leading}
\frac{k\partial_z\phi}{\phi}=-2c,\qquad
\frac{\partial_\tau\phi}{\phi}=10t_1+a-c^2,\qquad
k\partial_z a=b.
\end{equation}
\end{enumerate}
\end{proposition}
\begin{proof}
The Siegel operator maps $E_4^{(2)}$ and $E_6^{(2)}$ to $E_4$
and $E_6$. The first Fourier--Jacobi coefficients of the two cusp
forms are nonzero scalar multiples of the specified Jacobi cusp
forms \cite{EZ,vdG}. Now use
$E_4=12k^4t_2$, $E_6=-216k^6t_3$, and
$\Delta_{\mathrm{ell}}=-k^{12}\Delta$ to obtain
\eqref{eq:Igusa-leading}. Since $\Delta h$ does not vanish on
$\mathcal H$, division gives the limit. Finally,
\[
 \R_z\Delta=0,\quad \R_\tau\Delta=12t_1\Delta,\quad
 \R_zh=-2ch,\quad \R_\tau h=(a-c^2-2t_1)h.
\]
Pull these identities back by $\mathsf t$ and use
$\dd\mathsf t(k\partial_z)=\R_z\circ\mathsf t$ and
$\dd\mathsf t(\partial_\tau)=\R_\tau\circ\mathsf t$.
\end{proof}

\subsection{The separating degeneration}


Let \(\mathcal C\to S\) be a semistable family of genus two curves. Suppose each fiber over a smooth divisor \(D\subset S\) consists of two elliptic curves meeting at a single node. Work locally on \(D\), label the two components \(X_1,X_2\), and take the preimage of the node on each component as its origin. Write \(\T_1\) for the moduli space of elliptic curves with origin and a symplectic de Rham frame whose first element is holomorphic.

\begin{proposition}\label{prop:enhanced-sep-limit}
\begin{enumerate}
\item[(1)] Monodromy around \(D\) acts trivially on \(H^1\) of a smooth fiber. The limiting polarized Hodge structure is \(H^1(X_1)\oplus H^1(X_2)\). On \(D\), the extended Gauss-Manin connection is the direct sum of the two elliptic connections.
\item[(2)] If the extended frame restricts to a pair of such elliptic frames, the family defines a map \(D\to\T_1\times\T_1\) after any additional Weierstrass marking is forgotten. In homology bases chosen on the components, with both cycles and forms ordered by component, the boundary period matrix is \(P_2|_D=P_{X_1}\oplus P_{X_2}\).
\end{enumerate}
\end{proposition}

\begin{proof}
The separating loop bounds a subsurface of a smooth fiber. Its homology class is zero, so the Picard-Lefschetz formula gives trivial monodromy. For a nodal curve of compact type, specialization identifies the limiting polarized Hodge structure with the direct sum of the cohomologies of its components \cite[Section~5.4]{Hain}. This identification is compatible with parallel transport on \(D\), giving the assertion on connections. The two component frames determine the stated map to \(\T_1\times\T_1\). A class from either component has zero periods on cycles from the other, so the off diagonal blocks of the period matrix vanish.
\end{proof}

In the marked Siegel upper half space, the product of the two elliptic upper half planes is the locus \(\epsilon=\Omega_{12}=0\). Formula~\eqref{eq:separating-period-matrix} shows that \(\epsilon\) vanishes to first order when the two fixed pointed curves are plumbed. Witt's restriction formula \cite{Witt} and its leading coefficient version \cite[Section~9]{vdG} describe the Taylor expansion normal to this locus.

\begin{proposition}[Witt]\label{prop:separating-coefficients}
Let $\mathfrak F$ be a scalar holomorphic Siegel modular form of
integral weight $r$ for $\Sp_4(\Z)$, with trivial character. Write
\[
 \mathfrak F\begin{pmatrix}\tau_1&\epsilon\\
                         \epsilon&\tau_2\end{pmatrix}
       =\sum_{j\ge0}\psi_j(\tau_1,\tau_2)\epsilon^j.
\]
For \(n\ge0\), if \(\psi_j=0\) for all \(j<n\), then \(\psi_n\) is an elliptic
modular form of weight \(r+n\) in each variable. For every \(j\ge0\),
\(\psi_j=0\) unless \(j\equiv r\pmod 2\).
\end{proposition}
\begin{proof}
For \(\gamma_i=\left(\begin{smallmatrix}a_i&b_i\\c_i&d_i\end{smallmatrix}\right)
\in\SL_2(\Z)\), put \(j_i=c_i\tau_i+d_i\) and
\(J=j_1j_2-c_1c_2\epsilon^2\). Under the embedded pair
\((\gamma_1,\gamma_2)\), the new off diagonal entry is \(\epsilon/J\),
the new diagonal entries are \(\gamma_i\tau_i+O(\epsilon^2)\),
and the automorphy factor is \(J^r\). Since all lower coefficients vanish,
comparison of the terms of degree \(n\) gives
\[
 \psi_n(\gamma_1\tau_1,\gamma_2\tau_2)
   =j_1^{r+n}j_2^{r+n}\psi_n(\tau_1,\tau_2).
\]
In the Fourier expansion of \(\mathfrak F\), the indexing matrices are positive semidefinite. Their diagonal entries are nonnegative, so the expansion of \(\psi_n\) has no negative powers of either elliptic
Fourier parameter. Finally, the pair \((I,-I)\) fixes \(\tau_1,\tau_2\),
sends \(\epsilon\) to \(-\epsilon\), and has automorphy factor \((-1)^r\).
Comparison of coefficients gives the parity assertion.
\end{proof}

\bigskip
{\bf Acknowledgements.}
This project originated from discussions with Hossein Movasati about further  extending our previous work on quasi Jacobi forms of index zero \cite{CMVL},  which began two years ago. However, for a long time we were troubled by some  incompatibility in the computation of period matrices and the Gauss-Manin connection. During that period, we proposed several potential approaches;  see \cite[Chapter 14]{MovBook}. Later, Movasati gave up on this project for  some reason, and the author decided to continue the work and used the method of tangential regularization to  resolve the incompatibility.  I thank him for many helpful discussions and for his encouragement. I am also grateful to Gerard van der Geer for numerous helpful conversations.

The work was partially supported by the following grants: International Exchange and Development Program for Young Faculty at USTB, QNXM20250045 (Grant No.  06320202) the University of Science and Technology Beijing Foundation (Grant No. 00007886) and the Fundamental Research Funds for the Central Universities, FRF-BRA-26-007 (Grant No. 06108368).

\end{document}